\documentclass[10pt]{article}
\usepackage{amssymb}
\usepackage{amsthm}
\usepackage{graphicx}
\usepackage{graphics}
\usepackage{amsthm}

\input pictex.tex
\newcommand{\esp}{\hspace{0.01cm}}

\newcommand{\clo}{\mathrm{S}^1}

\theoremstyle{definition}

\newtheorem{thm}{Theorem}[section] 

\newtheorem{prop}[thm]{Proposition}

\newtheorem{lem}[thm]{Lemma}

\newtheorem{rem}[thm]{Remark}

\newtheorem{ex}[thm]{Example}

\input epsf

\begin{document}

\date{}
\author{Andr\'es Navas}

\title{Growth of groups and diffeomorphisms of the interval}
\maketitle

\vspace{-0.4cm}

\noindent{\bf Abstract.} We prove that, for all $\alpha \!>\! 0$, every finitely 
generated group of $C^{1+\alpha}$ diffeomorphisms of the interval with 
sub-exponential growth is almost nilpotent. Consequently, there is 
no group of $C^{1+\alpha}$ interval diffeomorphisms having intermediate 
growth. In addition, we show that the $C^{1+\alpha}$ regularity 
hypothesis for this assertion is essential by giving a $C^1$ 
counter-example.

\vspace{0.15cm}

\noindent{\bf MSC-class:} 20B27, 37E10, 37C85.

\vspace{0.54cm}

\noindent{\Large {\bf Introduction}}

\vspace{-0.15cm}

\vspace{0.5cm} A theory for groups of diffeomorphisms of the 
interval has been extensively developed by many authors (see for example 
\cite{Ko,Na-subexp,Na-solv,PT,serg,sz,Th,Ts,Ts-annals,tts,yoccoz}). One of the most 
interesting topics of this theory is the interplay between the differentiability 
class of the diffeomorphisms and the algebraic (as well as dynamical) properties 
of the group (and the action). For instance, as a consequence of the classical 
Bounded Distortion Principle, groups of $C^2$ diffeomorphisms appear to have 
a very rigid behavior. This is no longer true for subgroups of 
$\mathrm{Diff}_+^1([0,1])$, as it is well illustrated in the 
literature \cite{FF,Pix,Ts}. Major progress has recently been 
made in the understanding of the lost of rigidity in intermediate 
differentiability classes between $C^2$ and $C^1$ (see \cite{DKN}). 
The object of this work is to study the latter phenomenon for a 
remarkable class of groups, first introduced by R. Grigorchuk.

Given a finitely generated group (provided with a finite and symmetric system of 
generators), the growth function  assigns to each positive integer $n$ the number 
of elements of the group that may be written as a product of no more 
than $n$ generators. One says that the group has polynomial, exponential 
or intermediate growth, if its growth function has the 
corresponding asymptotic behaviour. (These notions do not depend on the choice of the 
finite system of generators.) A celebrated theorem by M. Gromov establishes that a 
group has polynomial growth if and only if it is almost nilpotent, {\em i.e.} if 
it contains a finite index nilpotent subgroup (see \cite{gr} and references therein). 
Typical examples of groups with exponential growth are those that contain free
semi-groups 
on two generators. (However, there exist groups with exponential growth and no free 
semi-group on two generators; see \cite{olsha}.) The difficult question (raised 
by J. Milnor \cite{milnor}) concerning the existence of groups with intermediate 
growth was positively answered by R. Grigorchuk in \cite{grigor-deg} (see also 
\cite{grigor2}). Some years later, one of his examples was realized (by R. Grigorchuk 
himself and A. Maki \cite{GrM}) as a subgroup of $\mathrm{Homeo}_+([0,1])$. 
The problem of improving the regularity for this embedding is at the core 
of this work. In the first part of this article we prove the following.

\vspace{0.4cm}

\noindent{\bf Theorem A.} {\em There exists a finitely generated subgroup of 
$\mathrm{Diff}_+^1([0,1])$ with intermediate growth.} 

\vspace{0.4cm}

This theorem solves by the negative a conjecture of \cite{GrM}. In fact, 
the group we consider turns out to be isomorphic to the group introduced by 
R. Grigorchuk in \cite{grigor2} and studied in more detail in \cite{GrM}. In the 
sequel, we will denote this group by $H$. We will prove more generally that, for 
every 
$C^1$-neighborhood $V$ of the identity map of $[0,1]$, there exists an embedding 
$H \hookrightarrow \mathrm{Diff}_+^1([0,1])$ sending some canonical system of 
generators of $H$ into $V$. This last issue is interesting because 
it is known for instance 
that subgroups of $\mathrm{Diff}_+^2(\clo)$ generated by elements 
near the identity (with respect to the $C^2$-topology) have very restrictive 
dynamical properties \cite{Na-duminy}.

The proof of Theorem A has two main technical ingredients. One is that, 
instead of embedding directly $H$ into $\mathrm{Diff}_+^1([0,1])$ 
(which seems to be very difficult), we construct a coherent sequence 
of embeddings of some almost nilpotent groups $H/H_n$ which in some 
sense converge to $H$. (The group $H$ turns out to be residually 
almost nilpotent.) An equicontinuity argument allows us to obtain, 
at the limit, the desired embedding. However, in order to apply 
this argument, it is necessary to ensure a uniform control for the 
derivatives of the generators of each group in the afore mentioned
sequence with respect to some fixed modulus of continuity. To do 
this, the other ingredient of the construction is to use a technique inspired 
by Chapter X of M. Herman's thesis \cite{Her}. This is related to the classical 
construction of $C^{1+\alpha}$ Denjoy counter-examples which, as explained 
to the author by F. Sergeraert, seems to go back to J. Milnor. 

The preceding method of proof is quite natural because 
the dynamics of (the canonical action of) the group $H$ has 
{\em infinitely many levels} (in the sense of \cite{hector2}; see also 
\cite[Section 8.3]{CClibro}), and a certain amount of regularity 
is lost when passing 
from one level to another ({\em i.e.} from the embedding of $H/H_n$ to that of 
$H/H_{n+1}$). In particular, at the limit we do not obtain an inclusion of $H$ into 
$\mathrm{Diff}_+^{1+\alpha}([0,1])$ for any $\alpha > 0$. (We get an uniform control 
for the derivatives of the generators only with respect to a logarithmic modulus of 
continuity.) And indeed, this issue is impossible, because of the following 
theorem.

\vspace{0.4cm}

\noindent{\bf Theorem B.} {\em For all $\alpha > 0$ every finitely generated 
subgroup of $\mathrm{Diff}_+^{1+\alpha}([0,1])$ with sub-exponential growth 
is almost nilpotent.}

\vspace{0.4cm}

This result is proved in the second part of this article (which is essentially 
independent of the first part), and holds more generally for (finitely generated) 
groups without free semi-groups on two generators. The proof relies on the rigidity 
theory for centralizers of diffeomorphisms of the interval. The foundations of this 
theory are related to the so-called Kopell Lemma \cite{Ko} (each Abelian subgroup 
of $\mathrm{Diff}_+^2([0,1[)$ either acts freely on $]0,1[$ or has a global 
fixed point therein), and Szekeres' theorem \cite{yoccoz} (the centralizer 
in $\mathrm{Diff}_+^1([0,1[)$ of every element of $\mathrm{Diff}^2_+([0,1[)$ 
without fixed points in $]0,1[$ is conjugate to the group of translations). 
For other classes of groups there is Plante-Thurston Theorem \cite{PT} 
(nilpotent subgroups of $\mathrm{Diff}_+^2([0,1[)$ are Abelian), and 
the classification of solvable subgroups of $\mathrm{Diff}_+^2([0,1[)$ 
obtained by the author in \cite{Na-solv} (see also \cite{Na-subexp}). 
Remark, however, that in the afore mentioned results, a $C^2$ regularity 
hypothesis is always assumed. (Or at least it is supposed that the maps are $C^1$ 
with derivatives having finite total variation.) The possibility of obtaining a 
result like Theorem B in intermediate regularity class was first suggested in 
\cite{DKN}. Let us mention that it is relatively simple to adapt the methods 
of proof to show that finitely generated subgroups 
of $\mathrm{Diff}_+^{1+\alpha}(\mathbb{R})$ or 
$\mathrm{Diff}^{1+\alpha}_+(\mathrm{S}^1)$ with 
sub-exponential growth (or without free semi-groups 
on two generators) are also almost nilpotent. 

Although no non Abelian nilpotent group can be contained in 
$\mathrm{Diff}_+^2([0,1[)$, a result of \cite{FF} establishes that every 
finitely generated torsion free nilpotent group can be seen as a subgroup of 
$\mathrm{Diff}_+^1([0,1])$. Using the methods of \cite{Ts}, it seems that the 
regularity of these inclusions can be improved up to the class $C^{1+\alpha}$ 
for every $\alpha < 1/(k-1)$, where $k$ is the nilpotence degree of the 
corresponding group. However, as a consequence of one of the main results 
of \cite{DKN}, these actions cannot be made $C^{1+\alpha}$ 
for any $\alpha > 1/(k-1)$. (And the same should be true for $\alpha = 1/(k-1)$.) 
Quite surprisingly, the proof of Theorem B is somehow different than the 
proof of this last statement. On the one hand it does not use the probabilistic 
techniques introduced in \cite{DKN} (see also \cite{tao}) 
to get control of distortion estimates 
in sharp intermediate differentiability classes. However, since the 
topological dynamics for the action is not prescribed {\em a priori}, 
it needs of an accurate study of the combinatorial properties for 
continuous actions of subexponential growth groups on the interval.

As it was previously recalled, among finitely generated groups the almost nilpotent 
ones are exactly those whose growth is polynomial \cite{gr}. Hence, Theorem B can 
be restated by saying that finitely generated sub-exponential growth subgroups 
of $\mathrm{Diff}_+^{1+\alpha}([0,1])$ have polynomial growth. As a direct 
corollary we obtain the following result, which implies that the 
conjecture of \cite{GrM} was true up to some $\alpha > 0$ !
 
\vspace{0.4cm}

\noindent{\bf Corollary.} {\em For all $\alpha > 0$ there is no finitely generated 
subgroup of $\mathrm{Diff}_+^{1+\alpha}([0,1])$ having intermediate growth.}

\vspace{0.4cm}

We conclude this Introduction with some remarks. According to the comments 
after Proposition 5.15 of \cite{DKN}, Grigorchuk-Maki's group $H$ seems to 
be the first example of a group of $C^1$ diffeomorphisms of the interval 
which cannot be seen as a group of $C^{1+\alpha}$ interval diffeomorphisms. 
On the other hand, it is well known that intermediate growth groups cannot 
appear as subgroups of Lie groups. These facts show that 
$\mathrm{Diff}_{+}^{1+\alpha}([0,1])$ is more appropriate 
than $\mathrm{Diff}^1([0,1])$ as an infinite dimensional model 
of a Lie group. Finally, note that the statements of Theorems A 
and B suggest a certain relationship with the classical Pesin Theory 
for diffeomorphisms with hyperbolic properties 
(see for instance \cite{katok}). Although these theorems 
seem to be of a different nature (our growth condition concerns the group 
and not the dynamics, and the control of distortion in class $C^{1+\alpha}$ is 
not related to any hyperbolicity of the action), possible further developments 
of a Pesin like theory for group actions could lead to a nice framework where 
these results appear as natural pieces.

\vspace{0.5cm}

\noindent{\bf Acknowledgments.} The author is indebted to L. Bartholdi, 
E. Breuillard, S. Crovisier, B. Deroin, A. Erschler, 
\'E. Ghys, T. Tsuboi and J. C. Yoccoz for fruitful discussions and 
their interest on the subject, and to J. Kiwi and J. Rivera-Letelier 
for several and valuable remarks and corrections to this article. 
Part of this work was supported by U. of Chile's DI-REIN Grant 06-01. 

%%%%%%%%%%%%%%%%%%%%%%%%%%%%%%%%%%%%%%%%%%%%%%%%%%%%%%%%%%%%%%%%%%%%%%%%%%%%%%%%%%%%%%%%%%%%%%%%%%
%%%%%%%%%%%%%%%%%%%%%%%%%%%%%%%%%%%%%%%%%%%%%%%%%%%%%%%%%%%%%%%%%%%%%%%%%%%%%%%%%%%%%%%%%%%%%%%%%%

\section{A group of $C^1$ interval diffeomorphisms with intermediate growth}

\hspace{0.35cm} The first half of this article is devoted to the proof of Theorem A. 
In \S \ref{cantor}, we recall the definition of the first Grigorchuk's group $G$ as well 
as Grigorchuk-Maki's group $H$ mainly as groups acting on spaces of sequences. In \S 
\ref{equiv}, we introduce several procedures for obtaining natural actions of $H$ by 
interval homeomorphisms. In \S \ref{modulo}, we study in detail the analytic 
properties for one of these procedures for getting some technical but quite 
useful estimates. Finally, in \S \ref{embedding}, we use these analytic 
estimates to obtain a faithful action of $H$ by $C^1$ 
diffeomorphisms of $[0,1]$. 
 
%%%%%%%%%%%%%%%%%%%%%%%%%%%%%%%%%%%%%%%%%%%%%%%%%%%%%%%%%%%%%%%%%%%%%%%%%%%%%%%%%%%%%%%%%%%%%%%%%%
%%%%%%%%%%%%%%%%%%%%%%%%%%%%%%%%%%%%%%%%%%%%%%%%%%%%%%%%%%%%%%%%%%%%%%%%%%%%%%%%%%%%%%%%%%%%%%%%%%

\subsection{Continuous actions on the interval and the Cantor set} 
\label{cantor}

\hspace{0.35cm} First Grigorchuk's group $G$ can be seen in many different ways: as 
the group generated by a finite automaton, as a group 
acting on the binary rooted tree $\mathcal{T}_2$, and as a group acting isometrically 
on the Cantor set $\{0,1\}^{\mathbb{N}}$. The last two points of view are essentially 
the same, since the boundary at infinity of $\mathcal{T}_2$ can be
identified with $\{0,1\}^{\mathbb{N}}$. 
Using the convention $(l_1,(l_2,l_3,\ldots)) = (l_1,l_2,l_3,\ldots)$ for $l_i \! \in
\! \{0,1\}$, the generators of $G$ are the elements
$\bar{\tt{a}},\bar{\tt{b}},\bar{\tt{c}},\bar{\tt{d}}$ 
whose actions on sequences 
$(l_1,l_2,l_3,\ldots)$ in $\{0,1\}^{\mathbb{N}}$ are defined recursively by 
$$\bar{\tt{a}}(l_1,l_2,l_3,\ldots) = (1-l_1,l_2,l_3,\ldots),$$ 
$$\bar{\tt{b}}(l_1,l_2,l_3,\ldots) = \left \{ \begin{array} {l} 
(l_1, \bar{\tt{a}}(l_2,l_3, \ldots )), \hspace{0.1cm} l_1 = 0,\\ 
(l_1, \bar{\tt{c}}(l_2,l_3,\ldots)), \hspace{0.1cm} l_1 = 1, \end{array} \right.$$
$$\bar{\tt{c}}(l_1,l_2,l_3,\ldots) = \left \{ \begin{array} {l} 
(l_1, \bar{\tt{a}}(l_2,l_3, \ldots )), \hspace{0.1cm} l_1 = 0,\\ 
(l_1, \bar{\tt{d}}(l_2,l_3,\ldots)), \hspace{0.1cm} l_1 = 1, \end{array} \right.$$
$$\bar{\tt{d}}(l_1,l_2,l_3,\ldots) = \left \{ \begin{array} {l} 
(l_1, l_2,l_3, \ldots ), \hspace{0.53cm} l_1 = 0,\\ 
(l_1, \bar{\tt{b}}(l_2,l_3,\ldots)), \hspace{0.1cm} l_1 = 1. \end{array} \right.$$
The action on $\mathcal{T}_2$ of the element $\bar{\tt{a}} \in G$ consists
of permuting the first two edges (and consequently, the trees rooted
on the final vertex of each one of those edges). Elements $\bar{\tt{b}}$, 
$\bar{\tt{c}}$ and $\bar{\tt{d}}$ fix the first two edges of $\mathcal{T}_2$, 
and their action on higher levels is illustrated in Figure 1 below.

\vspace{0.35cm}

%%%%%%%%%%%%%%%%%%%%%%%%%%%%%%%%%%%%%%%%%%%%%%%%%%%%%%%%%%%%%%%%%%%%%%%%%%%%%%%%%%%%%%%%%%%%%%%%%
%%%%%%%%%%%%%%%%%%%%%%%%%%%%%%%%%%%%%%%%%%%%%%%%%%%%%%%%%%%%%%%%%%%%%%%%%%%%%%%%%%%%%%%%%%%%%%%%%

\beginpicture

\setcoordinatesystem units <1cm,1cm>

%%%%%%%%%%%%%%%%%%%%%%%%%%%%%%%%%%%%%%%%%%%%%%%   a  
%%%%%%%%%%%%%%%%%%%%%%%%%%%%%%%%%%%%%%%%%%%%%%%

\plot 
0 0 
0.7 -0.8 /
\plot 
0 0 
-0.7 -0.8 /
\plot 
0.7 -0.8 
1 -1.5 /
\plot 
0.7 -0.8 
0.4 -1.5 /
\plot 
-0.7 -0.8 
-0.4 -1.5 /
\plot 
-0.7 -0.8 
-1 -1.5 /
\plot 
1 -1.5 
1.15 -1.9 /
\plot 
1 -1.5 
0.85 -1.9 /
\plot 
0.4 -1.5 
0.55 -1.9 /
\plot 
0.4 -1.5 
0.25 -1.9 /
\plot 
-0.4 -1.5 
-0.25 -1.9 /
\plot 
-0.4 -1.5 
-0.55 -1.9 /
\plot 
-1 -1.5 
-1.15 -1.9 /
\plot 
-1 -1.5 
-0.85 -1.9 /
\put{........................} at 0 -2.3

%%%%%%%%%%%%%%%%%%%%%%%%%%%%%%%%%%%%%%%%%%%%%%   b  
%%%%%%%%%%%%%%%%%%%%%%%%%%%%%%%%%%%%%%%%%%%%%%%%

\plot 
3.5 0 
4.2 -0.8 /
\plot 
3.5 0 
2.8 -0.8 /
\plot 
4.2 -0.8 
4.5 -1.5 /
\plot 
4.2 -0.8 
3.9 -1.5 /
\plot 
2.8 -0.8 
3.1 -1.5 /
\plot 
2.8 -0.8 
2.5 -1.5 /
\plot 
4.5 -1.5 
4.65 -1.9 /
\plot 
4.5 -1.5 
4.35 -1.9 /
\plot 
3.9 -1.5 
4.05 -1.9 /
\plot 
3.9 -1.5 
3.75 -1.9 /
\plot 
3.1 -1.5 
3.25 -1.9 /
\plot 
3.1 -1.5 
2.95 -1.9 /
\plot 
2.5 -1.5 
2.35 -1.9 /
\plot 
2.5 -1.5 
2.65 -1.9 /
\put{........................} at 3.5 -2.3

%%%%%%%%%%%%%%%%%%%%%%%%%%%%%%%%%%%%%%%%%%%%%%%%%%%%%%%%%%%%%%%%%%%%%%%%%%%%%%%%%%%%%%%%%%%%%%%%%%%

\plot 
0 0 
0.7 -0.8 /
\plot 
0 0 
-0.7 -0.8 /
\plot 
0.7 -0.8 
1 -1.5 /
\plot 
0.7 -0.8 
0.4 -1.5 /
\plot 
-0.7 -0.8 
-0.4 -1.5 /
\plot 
-0.7 -0.8 
-1 -1.5 /
\plot 
1 -1.5 
1.15 -1.9 /
\plot 
1 -1.5 
0.85 -1.9 /
\plot 
0.4 -1.5 
0.55 -1.9 /
\plot 
0.4 -1.5 
0.25 -1.9 /
\plot 
-0.4 -1.5 
-0.25 -1.9 /
\plot 
-0.4 -1.5 
-0.55 -1.9 /
\plot 
-1 -1.5 
-1.15 -1.9 /
\plot 
-1 -1.5 
-0.85 -1.9 /
\put{........................} at 0 -2.3

%%%%%%%%%%%%%%%%%%%%%%%%%%%%%%%%%%%%%%%%%%%%%%%   c  
%%%%%%%%%%%%%%%%%%%%%%%%%%%%%%%%%%%%%%%%%%%%%

\plot 
7 0 
7.7 -0.8 /
\plot 
7 0 
6.3 -0.8 /
\plot 
7.7 -0.8 
8 -1.5 /
\plot 
7.7 -0.8 
7.4 -1.5 /
\plot 
6.3 -0.8 
6.6 -1.5 /
\plot 
6.3 -0.8 
6 -1.5 /
\plot 
8 -1.5 
8.15 -1.9 /
\plot 
8 -1.5 
7.85 -1.9 /
\plot 
7.4 -1.5 
7.55 -1.9 /
\plot 
7.4 -1.5 
7.25 -1.9 /
\plot 
6.6 -1.5 
6.45 -1.9 /
\plot 
6.6 -1.5 
6.75 -1.9 /
\plot 
6 -1.5 
5.85 -1.9 /
\plot 
6 -1.5 
6.15 -1.9 /
\put{........................} at 7 -2.3
%%%%%%%%%%%%%%%%%%%%%%%%%%%%%%%%%%%%%%%%%%%%%%%   d  
%%%%%%%%%%%%%%%%%%%%%%%%%%%%%%%%%%%%%%%%%%%%%

\plot 
10.5 0 
11.2 -0.8 /
\plot 
10.5 0 
9.8 -0.8 /
\plot 
11.2 -0.8 
11.5 -1.5 /
\plot 
11.2 -0.8 
10.9 -1.5 /
\plot 
9.8 -0.8 
10.1 -1.5 /
\plot 
9.8 -0.8 
9.5 -1.5 /
\plot 
11.5 -1.5 
11.65 -1.9 /
\plot 
11.5 -1.5 
11.35 -1.9 /
\plot 
10.9 -1.5 
11.05 -1.9 /
\plot 
10.9 -1.5 
10.75 -1.9 /
\plot 
10.1 -1.5 
10.25 -1.9 /
\plot 
10.1 -1.5 
9.95 -1.9 /
\plot 
9.5 -1.5 
9.35 -1.9 /
\plot 
9.5 -1.5 
9.65 -1.9 /
\put{........................} at 10.5 -2.3

%%%%%%%%%%%%%%%%%%%%%%%%%%%%%%%%%%%%%%%%%%%%%%%%%%%%%%%%%%%%%%%%%%%%%%%%%%%%%%%%%%%%%%%%%%%%%%%%%

\put{$\bar{\tt{a}}$} at 0 0.6 
\put{$\curvearrowright$} at 0 0.2 
\put{$\longleftrightarrow$} at 0 -0.7 

\put{$\bar{{\tt b}}$} at 3.5 0.6 
\put{$\curvearrowright$} at 3.5 0.2 
\put{$\bar{\tt{a}}$} at 2.8 -0.2 
\put{$\bar{{\tt c}}$} at 4.3 -0.2 
\put{$\curvearrowright$} at 2.8 -0.5 
\put{$\curvearrowright$} at 4.3 -0.5 

\put{$\bar{{\tt c}}$} at 7 0.6 
\put{$\curvearrowright$} at 7 0.2 
\put{$\bar{\tt{a}}$} at 6.3 -0.2 
\put{$\bar{{\tt d}}$} at 7.8 -0.2 
\put{$\curvearrowright$} at 6.3 -0.5 
\put{$\curvearrowright$} at 7.8 -0.5 

\put{$\bar{{\tt d}}$} at 10.5 0.6 
\put{$\curvearrowright$} at 10.5 0.2 
\put{$id$} at 9.8 -0.2 
\put{$\bar{{\tt b}}$} at 11.3 -0.2 
\put{$\curvearrowright$} at 9.8 -0.5 
\put{$\curvearrowright$} at 11.3 -0.5 

\put{Figure 1} at 5.2 -2.8

\small
\put{$\bullet$} at 0 0 
\put{$\bullet$} at 3.5 0 
\put{$\bullet$} at 7 0 
\put{$\bullet$} at 10.5 0 
\put{$\bullet$} at 0.7 -0.8 
\put{$\bullet$} at 4.2 -0.8 
\put{$\bullet$} at 7.7 -0.8 
\put{$\bullet$} at 11.2 -0.8 
\put{$\bullet$} at -0.7 -0.8 
\put{$\bullet$} at 2.8 -0.8 
\put{$\bullet$} at 6.3 -0.8 
\put{$\bullet$} at 9.8 -0.8 

\put{$\bullet$} at 1 -1.5 
\put{$\bullet$} at 4.5 -1.5 
\put{$\bullet$} at 1 -1.5
\put{$\bullet$} at 11.5 -1.5 
\put{$\bullet$} at 10.9 -1.5
\put{$\bullet$} at -0.4 -1.5 
\put{$\bullet$} at 3.1 -1.5
\put{$\bullet$} at 6.6 -1.5 
\put{$\bullet$} at 9.5 -1.5 

\put{$\bullet$} at 8 -1.5
\put{$\bullet$} at 0.4 -1.5 
\put{$\bullet$} at 3.9 -1.5 
\put{$\bullet$} at 7.4 -1.5 
\put{$\bullet$} at 10.1 -1.5 
\put{$\bullet$} at -1 -1.5 
\put{$\bullet$} at 2.5 -1.5
\put{$\bullet$} at 6 -1.5 
\put{$\bullet$} at 9.5 -1.5 

\put{} at -3 0 

\endpicture

%%%%%%%%%%%%%%%%%%%%%%%%%%%%%%%%%%%%%%%%%%%%%%%%%%%%%%%%%%%%%%%%%%%%%%%%%%%%%%%%%%%%%%%%%%%%%%%%%
%%%%%%%%%%%%%%%%%%%%%%%%%%%%%%%%%%%%%%%%%%%%%%%%%%%%%%%%%%%%%%%%%%%%%%%%%%%%%%%%%%%%%%%%%%%%%%%%%

\vspace{0.25cm}

It can be shown that $G$ is a torsion group: each element has order a power of $2$ 
(see \cite{grigor1} or \cite{harpe}). The first example of a torsion free group 
with intermediate growth was given in \cite{grigor2}. Geometrically, the idea 
consists of replacing $\mathcal{T}_2$ by a rooted tree having vertices of 
infinite (countable) degree. In other terms, we consider the group $H$ acting 
on the space $\Omega = \mathbb{Z}^{\mathbb{N}}$ which is generated by the 
elements $\tt{a,b,c}$ and $\tt{d}$ defined recursively by 
$${\tt a} (l_1,l_2,l_3, \ldots) = (1+l_1,l_2,l_3, \ldots),$$
$${\tt b}(l_1,l_2,l_3,\ldots) = \left \{ \begin{array} {l} 
(l_1, {\tt a}(l_2,l_3, \ldots )), \hspace{0.1cm} l_1 \mbox{ even},\\ 
(l_1, {\tt c}(l_2,l_3,\ldots)), \hspace{0.1cm} l_1 \mbox{ odd}, \end{array} \right.$$
$${\tt c}(l_1,l_2,l_3,\ldots) = \left \{ \begin{array} {l} 
(l_1, {\tt a}(l_2,l_3, \ldots )), \hspace{0.1cm} l_1 \mbox{ even},\\ 
(l_1, {\tt d}(l_2,l_3,\ldots)), \hspace{0.1cm} l_1 \mbox{ odd}, \end{array} \right.$$
$${\tt d}(l_1,l_2,l_3,\ldots) = \left \{ \begin{array} {l} 
(l_1, l_2,l_3, \ldots ), \hspace{0.53cm} l_1 \mbox{ even},\\ 
(l_1, {\tt b}(l_2,l_3,\ldots)), \hspace{0.1cm} l_1 \mbox{ odd}. \end{array} \right.$$

The group $H$ preserves the lexicographic order on $\Omega$. It is then a left
orderable group 
\cite{ghys}, and so it can be realized as a group of orientation preserving
homeomorphisms of 
the interval. These facts were first established (by an indirect method) in
\cite{GrM}.

We next give an elementary proof of the fact that $H$ can be realized as a group of
bi-Lipschitz 
homeomorphisms of $[0,1]$. 
Fix a sequence $(\ell_i)_{i \in \mathbb{Z}}$ of positive 
numbers such that $\sum \ell_i = 1$ and 
$$\max \left\{ \frac{\ell_{i+1}}{\ell_i}, \frac{\ell_i}{\ell_{i+1}} \right\} 
\leq M < \infty \quad \mbox { for all } \quad i \in \mathbb{Z}.$$
We let $I_i$ denote the interval $]\sum_{j<i}\ell_j,\sum_{j\leq i}\ell_j[$. 
Let $f: [0,1] \rightarrow [0,1]$ be the orientation preserving 
homeomorphism sending each interval $I_i$ onto $I_{i+1}$ affinely.
Let $g$ be the orientation preserving affine homeomorphism sending 
$[0,1]$ onto $I_0$, and let us denote by $\lambda = 1/\ell_0$ 
the (constant) value of its derivative. Consider the maps 
$A,B,C$ and $D$ defined recursively on a dense subset of 
$[0,1]$ by setting $A(x) = f(x)$ and, for $x \in I_i$,
$$B(x) = \left \{ \begin{array} {l} 
f^i g A g^{-1} f^{-i}(x), \hspace{0.1cm} i \mbox{ even},\\ 
f^i g C g^{-1} f^{-i}(x), \hspace{0.1cm} i \mbox{ odd}, \end{array} \right.$$
$$C(x) = \left \{ \begin{array} {l} 
f^i g A g^{-1} f^{-i}(x), \hspace{0.1cm} i \mbox{ even},\\ 
f^i g D g^{-1} f^{-i}(x), \hspace{0.1cm} i \mbox{ odd}, \end{array} \right.$$
$$D(x) = \left \{ \begin{array} {l} 
x, \hspace{2,27cm} i \mbox { even},\\ 
f^i g B g^{-1} f^{-i}(x), \hspace{0.1cm} i \mbox{ odd}. \end{array} \right.$$

We claim that $A, B, C$ and $D$ are bi-Lipschitz homeomorphisms with 
bi-Lipschitz constant bounded above by $M$. Indeed, this is clear for $A$. 
For $B$, $C$ and $D$, this fact can be 
easily verified by induction. For example, 
if $x \in I_i$ for an even integer $i$, then 
$$B'(x) = \frac{(f^i)'(gAg^{-1}f^{-i}(x))}{(f^i)'(f^{-i}(x))} \cdot 
\frac{g'(Ag^{-1}f^{-i}(x))}{g'(g^{-1}f^{-i}(x))} \cdot A'(g^{-1}f^{-i}(x)),$$
and since $g'|_{[0,1]} = \lambda$ and $(f^i)'|_{I_0} = \ell_i/\ell_0$, we 
obtain $B'(x) = A'(g^{-1}f^{-i}(x)) \leq M$. Therefore, the maps $A,B,C$ and 
$D$ extend to bi-Lipschitz homeomorphisms of the whole interval $[0,1]$, 
and it is geometrically clear that they generate a group isomorphic to $H$. 
Remark finally that the constant $M$ may be chosen so near to $1$ as we want.

\begin{rem}
The fact that $H$ can be realized as a group of bi-Lipschitz homeomorphisms is 
not surprising. Indeed, a simple argument using the harmonic measure allows to 
show that if $\Gamma$ is any finitely generated subgroup of $\mathrm{Homeo}_+(\clo)$ 
(resp. of  $\mathrm{Homeo}_+([0,1])$), then $\Gamma$ is topologically conjugate to 
a group of bi-Lipschitz homeomorphisms of the circle (resp. of the interval): see 
\cite[Theorem D]{DKN}. Nevertheless, the fact that the Lipschitz constant 
of the generators of $H$ can be taken so near to $1$ as desired is a particular 
property of the group $H$. This property seems to be shared by any other (finitely 
generated) group of homeomorphisms of the interval or the circle without free 
semi-groups on two generators, but it is easy to construct examples showing 
that it does not hold in general.
\end{rem}

The preceding idea is not appropriate for obtaining an embedding of $H$ into 
$\mathrm{Diff}_+^1([0,1])$. Indeed, the discontinuities for the derivative 
repeat at each level of the action of $H$. In the subsequent section, we 
will give a method of construction to obtain such an embedding. For this 
we will have to renormalize suitably the geometry at each step. Denoting 
by $H_n$ the stabilizer of the level $n$ of the tree $\mathcal{T}_{\infty}$ 
for the action of $H$, we will construct embeddings of $H/H_n$ into 
$\mathrm{Diff}_+^1([0,1])$ in a coherent way and keeping some uniform control for 
the derivatives of generators; then using Arzel\'a-Ascoli Theorem, we will pass 
to the limit and obtain the desired embedding. Unfortunately, this method of 
construction will involve some technical issues. As a mater of fact, the action 
we will obtain is only semi-conjugate, but not conjugate, to the bi-Lipschitz 
action constructed above. (However, this seems to be a necessary condition 
for $C^1$-actions of $H$.)

\vspace{0.1cm}

\begin{rem} In \cite{grigor-deg}, R. Grigorchuk gives a general procedure for 
constructing groups of intermediate growth as groups acting on the dyadic rooted 
tree $\mathcal{T}_2$. It is not very difficult to see that the induced groups 
acting on the tree $\mathcal{T}_{\infty}$ by order preserving maps still have 
intermediate growth. For all of these induced groups, the methods of the first 
part of this work lead to realizations as subgroups of $\mathrm{Diff}_+^1([0,1])$.
\end{rem}
 
\begin{rem} A nice example of a group 
having {\em non uniform exponential growth} (that is, its 
exponential rate of growth is positive but becomes arbitrarily 
small under suitable changes of the system of generators) 
has been recently given by J. Wilson in \cite{wilson}. This 
group acts faithfully by automorphisms of a rooted tree. It 
would be interesting to know wether there exists an associated group 
of non uniform exponential growth acting on $\mathcal{T}_{\infty}$ 
by order preserving transformations. If this is the case, then 
the methods of the first part of this work should certainly 
lead to realizations as a subgroup of $\mathrm{Diff}_+^1([0,1])$.
\label{rem-wil}
\end{rem}

%%%%%%%%%%%%%%%%%%%%%%%%%%%%%%%%%%%%%%%%%%%%%%%%%%%%%%%%%%%%%%%%%%%%%%%%%%%%%%%%%%%%%%%%%%%%%%%%%%%
%%%%%%%%%%%%%%%%%%%%%%%%%%%%%%%%%%%%%%%%%%%%%%%%%%%%%%%%%%%%%%%%%%%%%%%%%%%%%%%%%%%%%%%%%%%%%%%%%%%
%%%%%%%%%%%%%%%%%%%%%%%%%%%%%%%%%%%%%%%%%%%%%%%%%%%%%%%%%%%%%%%%%%%%%%%%%%%%%%%%%%%%%%%%%%%%%%%%%%%

\subsection{Embeddings using equivariant families of homeomorphisms}
\label{equiv}

\hspace{0.35cm} Henceforth, we will deal only with {\em orientation 
preserving} homeomorphisms between intervals. A family 
$\{ \varphi_{u,v}: [0,u] \rightarrow [0,v]; u>0,v>0\}$ of such homeomorphisms 
will be called {\em equivariant} if for all $u>0,v>0$ and $w>0$, 
one has $\varphi_{v,w} \circ \varphi_{u,v} = \varphi_{u,w}$. Given 
such a family and two nondegenerate intervals $I=[x_1,x_2]$ and $J=[y_1,y_2]$, 
we let $\varphi(I,J):I \rightarrow J$ denote the homeomorphism defined by 
$$\varphi(I,J)(x) = \varphi_{x_2-x_1,y_2-y_1}(x-x_1) + y_1.$$
Remark that $\varphi(I,I)$ is forced to be the identity map.

The simplest family of equivariant homeomorphisms is the one consisting 
of the affine maps $\varphi_{u,v}(x)~=~vx/u$. However, this family is not
adequate if we want to fit maps smoothly together. Let us then introduce 
a general and simple procedure for constructing families of equivariant 
homeomorphisms as follows. Let $\{ \varphi_{u}: \mathbb{R} \rightarrow
]0,u[; u>0 \}$ be any family of homeomorphisms. 
Define $\varphi_{u,v}: ]0,u[ \rightarrow ]0,v[$ by 
$\varphi_{u,v} = \varphi_{v} \circ \varphi_u^{-1}$. We have
$$\varphi_{v,w} \circ \varphi_{u,v} = (\varphi_w \circ \varphi_v^{-1}) \circ 
(\varphi_v \circ \varphi_u^{-1}) = \varphi_w \circ \varphi_u^{-1} =
\varphi_{u,w}.$$  
Thus, extending $\varphi_{u,v}$ continuously to the whole 
interval $[0,u]$ by setting $\varphi_{u,v}(0) = 0$ and 
$\varphi_{u,v}(u)=v$, we obtain the desired equivariant family. 

\vspace{0.05cm}

\begin{ex} Let $\varphi_u: \mathbb{R} \rightarrow ]0,u[$ be given by  
$$\varphi_u(x) = \frac{1}{\pi} \int_{-\infty}^x 
\frac{ds}{s^2+(1/u)^2} = \frac{u}{2} + \frac{u}{\pi} \arctan(ux).$$
The corresponding equivariant family $\{ \varphi_{u,v} : [0,u] \rightarrow [0,v];
u>0,v>0 \}$ will be essential in what follows. This family 
was introduced by J. C. Yoccoz, and it has been already used in 
\cite{FF}. The regularity properties of the maps $\varphi_{u,v}$ 
will be studied in \S \ref{modulo}. 
\label{familia-de-yoccoz}
\end{ex}

\vspace{0.05cm}

Now fix any equivariant family of homeomorphisms 
$\{ \varphi_{u,v}: [0,u] \rightarrow [0,v]; u>0,v>0 \}$. 
For each $n \in \mathbb{N}$ and each $(l_1,\ldots,l_n) \in \mathbb{Z}^n$, let us 
consider a non degenerate closed interval 
$I_{l_1,\ldots,l_n}\!=\![a_{l_1,\ldots,l_n},b_{l_1,\ldots,l_n}]$ 
and a (perhaps degenerate) closed interval 
$J_{l_1,\ldots,l_n}\!=\![c_{l_1,\ldots,l_n},d_{l_1,\ldots,l_n}]$, 
both contained in some interval $[0,T]$. Let us suppose that the 
following conditions are satisfied (see Figure 2):

\noindent{(i) $\sum_{l_1 \in \mathbb{Z}} |I_{l_1}| = T$ (where $| \cdot |$ denotes
the length of the 
corresponding interval),}

\noindent{(ii) $a_{l_1,\ldots,l_n} < c_{l_1,\ldots,l_n} \leq d_{l_1,\ldots,l_n} =
b_{l_1,\ldots,l_n}$, 
so in particular $J_{l_1,\ldots,l_n} \subset I_{l_1,\ldots,l_n}$,}

\noindent{(iii) $b_{l_1,\ldots,l_{n-1},l_n} = a_{l_1,\ldots,l_{n-1},1+l_n}$,}

\noindent{(iv) $\lim_{l_n \rightarrow -\infty} a_{l_1,\ldots,l_{n-1},l_n} =
a_{l_1,\ldots,l_{n-1}}$,}

\noindent{(v) $\lim_{l_n \rightarrow \infty} a_{l_1,\ldots,l_{n-1},l_n} =
c_{l_1,\ldots,l_{n-1}}$,}

\noindent{(vi) $\lim_{n \rightarrow \infty} \sup_{(l_1,\ldots,l_n) \in \mathbb{Z}^n} 
|I_{l_1,\ldots,l_n}| = 0$.}

\vspace{0.5cm}

%%%%%%%%%%%%%%%%%%%%%%%%%%%%%%%%%%%%%%%%%%%%%%%%%%%%%%%%%%%%%%%%%%%%%%%%%%%%%%%%%%%%%%%%%%%%%%
\beginpicture
\setcoordinatesystem units <1cm,1cm>

\putrule from 0 0 to 12 0 
\putrule from 0 1.5 to 5 1.5
\putrule from 7 1.5 to 12 1.5 
\putrule from 5.3 0.7 to 8 0.7
\putrule from 9.6 0.7 to 12 0.7
\putrule from 2.15 -0.5 to 3.12 -0.5

\put{$J_{l_1,\ldots,l_n}$} at 8.8 0.7

\begin{small}
\put{$I_{l_1,\ldots,l_n,l_{n\!+\!1}}$} at 2.6 -1
\put{$a_{l_1\!,\!\ldots\!,\!l_n\!,\!l_{n\!+\!1}}$} at 1.5 0.5 
\put{$b_{l_1\!,\!\ldots\!,\!l_n\!,\!l_{n\!+\!1}}$} at 3.3 0.5
\end{small}

\begin{tiny}
\put{$|$} at 1.8 0
\put{$|$} at 2.8 0
\put{$|$} at 1.8 -0.5
\put{$|$} at 2.8 -0.5
\end{tiny}

\begin{Large}
\put{$a_{l_1,\ldots,l_n}$} at -0.8 -0.9
\put{$d_{l_1,\ldots,l_n} = b_{l_1,\ldots,l_n}$} at 11.6 -0.9 
\put{$|$} at -0.6 0 
\put{$|$} at 11.4 0 
\put{$c_{l_1,\ldots,l_n}$} at 4.6 -0.9
\put{$|$} at 4.65 0
\put{$I_{l_1,\ldots,l_n}$} at 5.45 1.5
\end{Large}

\put{Figure 2} at 5.9 -1.8

\put{} at -1.77 2.2

\endpicture

%%%%%%%%%%%%%%%%%%%%%%%%%%%%%%%%%%%%%%%%%%%%%%%%%%%%%%%%%%%%%%%%%%%%%%%%%%%%%%%%%%%%%%%%%%%%%%

\vspace{0.5cm}

\noindent{Note that}
\begin{equation}
|J_{l_1,\ldots,l_n}|+\sum_{l_{n+1} \in \mathbb{Z}}
|I_{l_1,\ldots,l_n,l_{n+1}}|=|I_{l_1,\ldots,l_n}|.
\label{defjota}
\end{equation}

For each $n \in \mathbb{N}$, we will define homeomorphisms $A_n, B_n, C_n$ and 
$D_n$ in such a way that the group generated by them will be isomorphic to
$H/H_n$. For this, let us consider the 
homomorphisms $\phi_0$ and $\phi_1$ from the subgroup of $H$ generated by 
$\tt b$, $\tt c$ and $\tt d$ into $H$ defined by 
$$\phi_0({\tt b}) = {\tt a}, \quad \phi_0({\tt c}) = {\tt a}, \quad \phi_0({\tt d})
= id, 
\qquad \mbox{ and } \qquad \phi_1({\tt b}) = {\tt c}, \quad \phi_1({\tt c}) = {\tt d}, 
\quad \phi_1({\tt d}) = {\tt b}.$$

\noindent{\underbar{Definition of $A_n$}}

\noindent{-- If $p \in J_{l_1,\ldots,l_i}$ for some $i < n$, let 
$A_n(p) = \varphi(J_{l_1,l_2,\ldots,l_i},J_{1+l_1,l_2,\ldots,l_i})(p)$.}

\noindent{-- If $p \in I_{l_1,\ldots,l_n}$, 
let $A_n(p) = \varphi(I_{l_1,l_2,\ldots,l_n},I_{1+l_1,l_2,\ldots,l_n})(p)$.}

\vspace{0.15cm}

\noindent{\underbar{Definition of $B_n$}}

Suppose that $p \! \in ]0,1[$ belongs to $I_{l_1,\ldots,l_n}$, 
and denote the corresponding sequence reduced modulo 2 by 
$(\bar{l}_1,\ldots,\bar{l}_n) \in \{0,1\}^n$.

\noindent{-- If $\phi_{\bar{l}_1}({\tt b})$, $\phi_{\bar{l}_2} \phi_{\bar{l}_1}({\tt
b})$, 
$\ldots$, $\phi_{\bar{l}_n} \cdots \phi_{\bar{l}_2} \phi_{\bar{l}_1}({\tt b})$ are 
well defined, let $B_n(p) = p$.}

\noindent{-- Otherwise, we denote the smallest integer $i \leq n$ such that 
$\phi_{\bar{l}_i} \cdots \phi_{\bar{l}_2} \phi_{\bar{l}_1}({\tt b})$ 
is not defined by $i(p)$.}

$\cdot$ If $p \in J_{l_1,\ldots,l_j}$ for some $j < i(p)$, let $B_n(p) = p$.

$\cdot$ If $p \in J_{l_1,\ldots,l_{i(p)},\ldots,l_j}$ for some $i(p) \leq j < n$, let 
$B_n(p) = 
\varphi(J_{l_1,\ldots,l_{i(p)},\ldots,l_j},J_{l_1,\ldots,1+l_{i(p)},\ldots,l_j})(p)$.

$\cdot$ If $p \in I_{l_1,\ldots,l_n}$, let $B_n(p) =
\varphi(I_{l_1,\ldots,l_{i(p)},\ldots,l_n},I_{l_1,\ldots,1+l_{i(p)},\ldots,l_n})(p)$.

\vspace{0.2cm}

The definitions of $C_n$ and $D_n$ are similar to that 
of $B_n$. Clearly, the maps $A_n,B_n,C_n$ and $D_n$ extend to 
homeomorphisms of $[0,T]$. The fact that they generate a group isomorphic to $H/H_n$
is 
geometrically clear and follows easily from the equivariant properties of the maps
$\varphi_{u,v}$. 
Moreover, condition (vi) implies that the sequences of maps $A_n$, $B_n$, $C_n$ and
$D_n$ converge 
to limit homeomorphisms $A,B,C$ and $D$ respectively, which generate a group
isomorphic to $H$.

\vspace{0.15cm}

\begin{ex} Given a sequence $(\ell_i)_{i \in \mathbb{Z}}$ of positive numbers such 
that $\sum \ell_i = 1$, define $|I_{l_1,\ldots,l_n}|$ and $|J_{l_1,\ldots,l_n}|$ by 
\esp \esp $|J_{l_1,\ldots,l_n}| = 0$ \esp \esp and \esp \esp  
$|I_{l_1,\ldots,l_n}| = \ell_{l_1} \cdots \ell_{l_n}$ 
respectively. If we carry 
out the preceding construction (for $T = 1$) using the equivariant family of affine 
maps $\varphi_{u,v}(x) = ux/v$, then we recover the embedding of $H$ into the 
group of bi-Lipschitz homeomorphisms of the interval constructed at the end 
of \S \ref{cantor} (under the assumptions that $\ell_{i+1}/\ell_i \leq M$ 
and $\ell_i/\ell_{i+1} \leq M$ for all $i \in \mathbb{Z}$). 
\label{emb-lip}
\end{ex}

%%%%%%%%%%%%%%%%%%%%%%%%%%%%%%%%%%%%%%%%%%%%%%%%%%%%%%%%%%%%%%%%%%%%%%%%%%%%%%%%%%%%%%%%%%%%%%%%%%%
%%%%%%%%%%%%%%%%%%%%%%%%%%%%%%%%%%%%%%%%%%%%%%%%%%%%%%%%%%%%%%%%%%%%%%%%%%%%%%%%%%%%%%%%%%%%%%%%%%%
%%%%%%%%%%%%%%%%%%%%%%%%%%%%%%%%%%%%%%%%%%%%%%%%%%%%%%%%%%%%%%%%%%%%%%%%%%%%%%%%%%%%%%%%%%%%%%%%%%%

\subsection{Modulus of continuity for the derivatives}
\label{modulo}

\hspace{0.35cm} Let $\sigma: [0,1] \rightarrow [0,\sigma(1)]$ be an increasing 
homeomorphism. A continuous map $\psi: [0,1] \rightarrow \mathbb{R}$ is
$\sigma$-continuous 
if there exists $M < \infty$ such that, for all $x \neq y$ in $[0,1]$, 
$$\left| \frac{\psi(x) - \psi(y)}{\sigma(|x-y|)} \right| \leq M.$$
We denote the supremum of the left hand side expression by $\|\psi\|_{\sigma}$, 
and we call it the $\sigma$-norm of $\psi$. The interest in the notion 
of $\sigma$-continuity relies on the obvious fact that, if $(\psi_n)$ 
is a sequence of functions defined on $[0,1]$ such that 
$$\sup_{n \in \mathbb{N}} \| \psi_n \|_{\sigma} < \infty,$$
then $(\psi_n)$ is an equicontinuous sequence.

\begin{ex} For $\sigma(s) = s^{\alpha}$, with $0\!<\!\alpha\!<\! 1$, the 
notions of $\sigma$-continuity and $\alpha$-H\"older continuity coincide. 
\label{wholder}
\end{ex}

\begin{ex} For $\varepsilon > 0$ suppose that $\sigma = \sigma_{\varepsilon}$ is such 
that $\sigma_{\varepsilon}(s) = s \log(1/s)^{1+\varepsilon}$ for $s$ small. If a map 
$\varphi$ is $\sigma_{\varepsilon}$-continuous, then it is 
$\alpha$-H\"older continuous for all 
$0 < \alpha < 1$. Indeed, it is easy to verify that 
$$s \left( \log \Big(\frac{1}{s} \Big) \right)^{1+\varepsilon} 
\leq C_{\varepsilon,\alpha} s^{\alpha}, \qquad \mbox{where} \quad
C_{\varepsilon,\alpha} 
= \frac{1}{e^{1+\varepsilon}} \left( \frac{1+\varepsilon}{1-\alpha}
\right)^{1+\varepsilon}.$$
We remark that the map $s \mapsto s \log(1/s)^{1+\varepsilon}$ is not 
Lipschitz. As a consequence, $\sigma_{\varepsilon}$-continuity for a 
function does not imply that the function is Lipschitz. 
\label{we}
\end{ex}

\begin{ex} A modulus of continuity $\sigma$ satisfying $\sigma(s) = 1/\log(1/s)$ 
for $s$ small enough is weaker than any H\"older modulus $s \mapsto s^{\alpha}$, with 
$\alpha > 0$. Nevertheless, such a modulus will be essential for our construction.
\label{wlog}
\end{ex}

We will now investigate several upper bounds with respect to some moduli 
of continuity for the derivatives of maps in Yoccoz's family (see Example 
\ref{familia-de-yoccoz}). Letting $y = \varphi_u^{-1}(x)$, we have 
$$\varphi_{u,v}'(x) = \varphi_v'(y)(\varphi_u^{-1})'(x) = 
\frac{\varphi_v'(y)}{\varphi_u'(y)} = \frac{y^2 + 1/u^2}{y^2 + 1/v^2}.$$
Note that, when $x \rightarrow 0$ (resp. $x \rightarrow u$), we have
that $y \rightarrow -\infty$ (resp. $y \rightarrow +\infty$), 
and $\varphi_{u,v}'(x) \rightarrow 1$. Therefore, the map 
$\varphi_{u,v}$ extends to a $C^1$ diffeomorphism 
from $[0,u]$ to $[0,v]$ which is tangent to the identity at the end points
of $[0,u]$. Moreover, for 
$u \geq v$ (resp. $u \leq v$), the function $s \mapsto \frac{s^2 + 1/u^2}{s^2 +
1/v^2}$ attains 
its minimum (resp. maximum) value at $s = 0$. Since this value is equal to
$v^2/u^2$, we have 
$$\sup_{x \in [0,u]} |\varphi_{u,v}'(x) - 1| = \left| \frac{v^2}{u^2} - 1 \right|.$$ 

For the second derivative of $\varphi_{u,v}$ we have 
\begin{small}
$$\varphi_{u,v}''(x) = \frac{d \varphi_{u,v}'(x)}{dy} \cdot \frac{dy}{dx} = 
\frac{2y(y^2+1/v^2) - 2y(y^2+1/u^2)}{(y^2+1/v^2)^2} \pi (y^2 + 1/u^2) = 
\pi \frac{y^2 + 1/u^2}{(y^2 + 1/v^2)^2} \left[ 2y \Big( \frac{1}{v^2} -
\frac{1}{u^2} \Big) \right].$$
\end{small}Therefore,
$$|\varphi_{u,v}''(x)| = \pi \frac{y^2 + 1/u^2}{y^2 + 1/v^2} \cdot 
\frac{|2y(1/v^2 - 1/u^2)|}{y^2 + 1/v^2}.$$ 
It follows from this equality that $\varphi_{u,v}$ is a $C^2$ diffeomorphism, 
with $\varphi_{u,v}''(0) \!=\! \varphi_{u,v}''(u)\!=\!0$. Moreover, the 
inequality $\frac{2|y|}{y^2 +t^2}\leq \frac{1}{t}$ applied to $t = 1/v$ yields 
$$\left| \varphi_{u,v}''(x) \right| \leq 
\pi \frac{y^2 + 1/u^2}{y^2 + 1/v^2} \left| \frac{1}{v^2} - \frac{1}{u^2} \right|v.$$ 
For $u \leq v$, this implies 
$$\left| \varphi_{u,v}''(x) \right| \leq 
\pi \frac{v^2}{u^2}\left( \frac{v^2 - u^2}{u^2v^2} \right)v = 
\frac{\pi v}{u^2} \left( \frac{v^2}{u^2}-1 \right).$$ 
So, if $u \leq v \leq 2u$, then 
$$\left| \varphi_{u,v}''(x) \right| \leq 6\pi \left| \frac{v}{u} - 1 \right|
\frac{1}{u}.$$ 
Analogously, if $2v \geq u \geq v$, then 
$$\left| \varphi_{u,v}''(x) \right| \leq 
\frac{\pi}{v} \left(1 - \frac{v^2}{u^2} \right) 
\leq 2\pi \left| \frac{v}{u}-1 \right| \frac{1}{v} 
\leq 4\pi \left| \frac{v}{u} - 1 \right| \frac{1}{u}.$$ 
Thus, in both cases, we have
\begin{equation}
\left| \varphi_{u,v}''(x) \right| \leq 6 \pi \left| \frac{v}{u} - 1 \right|
\frac{1}{u}.
\label{tres}
\end{equation} 
The last inequality and the next elementary proposition show that the family of 
maps $\varphi_{u,v}$ is in some sense optimal (at least among families of maps 
which are tangent to the identity at the end points; see Remark \ref{nt}).

\vspace{0.2cm}

\begin{prop} {\em If $\varphi: [0,u] \rightarrow [0,v]$ is a 
$C^2$ diffeomorphism such that $\varphi'(0)=\varphi'(u)=1$, then there exists a point 
$y \! \in ]0,u[$ such that} $$|\varphi''(y)|\geq \frac{2}{u} \left|
\frac{v}{u}-1\right|.$$
\end{prop}

\noindent{\bf Proof.} Let us suppose that $v \geq u$ (the case where $v \leq u$ 
is similar). Since $\varphi(0)=0$ and $\varphi(u)=v$, there exists some point 
$x \! \in ]0,u[$ such that $\varphi'(x) \geq v/u$. There are two cases:

\vspace{0.08cm}

\noindent{(i) $x \geq u/2$: the Mean Value Theorem gives some point $y \! \in ]x,u[$
such that}
$$|\varphi''(y)| = \frac{\varphi'(x)-\varphi'(u)}{u-x} \geq \frac{2}{u}|\varphi'(x)-1| 
\geq \frac{2}{u} \left| \frac{v}{u} - 1\right|;$$

\noindent{(ii) $x \leq u/2$: again, there exists $y \! \in ]0,x[$ such that} 
$$|\varphi''(y)| = \frac{\varphi'(x)-\varphi'(0)}{x} \geq \frac{2}{u}|\varphi'(x)-1| 
\geq \frac{2}{u} \left| \frac{v}{u} - 1\right|.$$

\vspace{0.28cm}

The next lemma should be compared to \S 3.17 of Chapter X of \cite{Her}. Note 
that the moduli from Examples \ref{wholder}, \ref{we} and \ref{wlog} can be taken 
satisfying the decreasing hypothesis on the function $s\mapsto \sigma(s)/s$ .

\vspace{0.15cm}

\begin{lem} {\em Suppose that $\sigma$ is a modulus of continuity such that 
the function $s\mapsto \sigma(s)/s$ is decreasing. If $u>0$ and $v>0$ 
satisfy $u/v \leq 2$, $v/u \leq 2$, and} 
$$\left| \frac{v}{u} - 1 \right| \frac{1}{\sigma(u)} \leq M,$$ 
{\em then the $\sigma$-norm of $\varphi_{u,v}'$ is less than or equal to $6\pi M$.}
\label{simple}
\end{lem}

\noindent{\bf Proof.} By inequality (\ref{tres}), for all $x \in [0,u]$ we have
$$|\varphi_{u,v}''(x)| \leq \frac{6 \pi M \sigma(u)}{u}.$$
If $y < z$ are points in $[0,u]$ then there exists  $x \in [y,z]$ such that 
$\varphi_{u,v}'(z) - \varphi'(y) = \varphi''_{u,v}(x) (z-y)$. Since 
$s \mapsto \sigma(s)/s$ is a decreasing function and $z-y \leq u$, this gives 
$$\left| \frac{\varphi_{u,v}'(z) - \varphi'_{u,v}(y)}{\sigma(z-y)} \right| =
|\varphi''_{u,v}(x)| \left| \frac{z-y}{\sigma(z-y)} \right| \leq 
|\varphi''_{u,v}(x)| \left| \frac{u}{\sigma(u)} \right| \leq 6 \pi M.$$
This finishes the proof of the lemma.

\vspace{0.15cm}

\begin{rem} In the case where 
non tangencies to the identity at the end points are allowed, there are 
equivariant procedures to construct maps with slightly better regularity properties
than 
those of the maps in Yoccoz's family. This is related to the famous Pixton's actions 
\cite{Pix}, for which an alternative and precious reference is \cite{Ts}. 
Nevertheless, for our construction we do not need a sequence of optimal 
embeddings of the groups $H / H_n$, but only a sequence of ``good enough'' embeddings 
which allow to preserve the differentiability when passing to the limit. For this
reason, 
we will not use the sharp constructions of \cite{Ts}, for which the computations are 
much more involved.
\label{nt}
\end{rem}

\vspace{0.1cm}

We finish this section with an elementary technical lemma which 
will be useful to control the $\sigma$-norm of the derivative of a map 
obtained by fitting together many diffeomorphisms defined on sub-intervals. 
The zero Lebesgue measure hypothesis below will be trivially satisfied 
in our constructions because the corresponding sets will be countable.

\vspace{0.15cm}

\begin{lem} {\em Let $\{I_n: n \in \mathbb{N}\}$ be a family of closed intervals in 
$[0,1]$ having disjoint interiors and such that the complement of their union 
has zero Lebesgue measure. Suppose 
that $\varphi$ is a homeomorphism of $[0,1]$ such that its restrictions to each 
interval $I_n$ are $C^{1+\sigma}$ diffeomorphisms which are $C^1$-tangent to the
identity 
at both end points of $I_n$ and whose derivatives have $\sigma$-norms bounded above 
by a constant $M$. Then $\varphi$ is a $C^{1+\sigma}$ diffeomorphism of the whole 
interval $[0,1]$, and the $\sigma$-norm of its derivative is less than or 
equal to $2M$.}
\label{pegar}
\end{lem}

\noindent{\bf Proof.} Let $x\!<\!y$ be two points of $\cup_{n \in \mathbb{N}} I_n$. 
If they belong to the same interval $I_n$ then, by hypothesis, 
$$\left| \frac{\varphi'(y) - \varphi'(x)}{\sigma(y-x)} \right| \leq M.$$
Suppose now that $x \in I_i = [x_i,y_i]$ and 
$y \in I_j = [x_j,y_j]$, with $y_i \leq x_j$. In this case,
\begin{eqnarray*}
\left| \frac{\varphi'(y)-\varphi'(x)}{\sigma(y-x)}
\right| 
&=& \left|\frac{(\varphi'(y)-1)+(1-\varphi'(x))}{\sigma(y-x)}\right| \\
&\leq& \left|\frac{\varphi'(y)-\varphi'(x_j)}{\sigma(y-x)} \right|
+ \left| \frac{\varphi'(y_i)-\varphi'(x)}{\sigma(y-x)}\right| \\ 
&\leq& M \left[ \frac{\sigma(y-x_j)}{\sigma(y-x)} +
\frac{\sigma(y_i-x)}{\sigma(y-x)} \right] \\
&\leq& 2M.
\end{eqnarray*}
The map $x\mapsto \varphi'(x)$ is then uniformly continuous on the dense set
$\cup_{n\in\mathbb{N}} 
I_n$, and so it extends to some continuous function on $[0,1]$ having $\sigma$-norm
bounded 
above by $2M$. Finally, since the complementary set of $\cup_{n\in\mathbb{N}} I_n$ 
has zero Lebesgue measure, the Fundamental Theorem of Calculus shows that this 
continuous function coincides (everywhere) with the derivative of $\varphi$. 

%%%%%%%%%%%%%%%%%%%%%%%%%%%%%%%%%%%%%%%%%%%%%%%%%%%%%%%%%%%%%%%%%%%%%%%%%%%%%%%%%%%%%%%%%%%%%%%%%%%
%%%%%%%%%%%%%%%%%%%%%%%%%%%%%%%%%%%%%%%%%%%%%%%%%%%%%%%%%%%%%%%%%%%%%%%%%%%%%%%%%%%%%%%%%%%%%%%%%%%
%%%%%%%%%%%%%%%%%%%%%%%%%%%%%%%%%%%%%%%%%%%%%%%%%%%%%%%%%%%%%%%%%%%%%%%%%%%%%%%%%%%%%%%%%%%%%%%%%%%

\subsection{The embedding of $\mathbf{\it H}$ into 
$\mathrm{Diff}_+^1([0,1])$}
\label{embedding}

\hspace{0.35cm} In this section, $\sigma$ will denote a fixed modulus of continuity
satisfying 
$\sigma(s) = 1/\log(1/s)$ for $s$ small enough (namely, for $s \leq 1/e$), and such
that the 
function $s \mapsto \sigma(s)/s$ is decreasing. We 
will prove that there exist embeddings $H \hookrightarrow
\mathrm{Diff}_+^{1+\sigma}([0,1])$ 
sending the generators $\tt a,b,c,d$ of $H$ (and their inverses) to diffeomorphisms 
so near as we want (in the $C^{1+\sigma}$-topology) to the identity map. To do that, 
fix any number $M>0$, and for each 
$k \in \mathbb{N}$ define $T_k = \sum_{i \in \mathbb{Z}} \frac{1}{(|i|+k)^2} <
\infty$. 
Consider an increasing sequence $(k_n)$ of positive integer numbers such that $k_1
\geq 4$. 
For $n \in \mathbb{N}$ and $(l_1,\ldots,l_n) \in \mathbb{Z}^n$ let
$$|I_{l_1,\ldots,l_n}| = \frac{1}{(|l_1|+ \cdots |l_n|+k_n)^{2n}}.$$
Note that 
\begin{eqnarray*}
\sum_{l_{n+1} \in \mathbb{Z}} |I_{l_1,\ldots,l_n,l_{n+1}}| 
&=& \sum_{l_{n+1}\in \mathbb{Z}} \frac{1}{(|l_1|+\cdots |l_n| + |l_{n+1}| +
k_{n+1})^{2n+2}} \\
&\leq& 2 \int_{|l_1| +\cdots +|l_n|+k_{n+1}-1}^{\infty} \frac{ds}{s^{2n+2}} \\
&=& \frac{2}{2n+1} \cdot \frac{1}{(|l_1|+\cdots +|l_n| + k_{n+1} - 1)^{2n+1}}.
\end{eqnarray*} 
Thus
\begin{small}
$$\frac{\sum_{l_{n+1} \in \mathbb{Z}}
|I_{l_1,\ldots,l_n,l_{n+1}}|}{|I_{l_1,\ldots,l_n}|} \leq 
\frac{2}{2n+1}\cdot \frac{(|l_1| + \cdots |l_n| +
k_n)^{2n}}{(|l_1|+\cdots+|l_n|+k_{n+1}-1)^{2n+1}} 
\leq \frac{2}{(2n+1)(|l_1| + \cdots + |l_n| + k_n)}.$$
\end{small}In 
particular, we can define $|J_{l_1,\ldots,l_n}|$ by (\ref{defjota}), that is 
$$|J_{l_1,\ldots,l_n}| = |I_{l_1,\ldots,l_n}| - \sum_{l_{n+1} \in \mathbb{Z}}
|I_{l_1,\ldots,l_n,l_{n+1}}|,$$
and for this choice we have 
\begin{equation}
|I_{l_1,\ldots,l_n}| \geq |J_{l_1,\ldots,l_n}| \geq 
\left( 1 - \frac{2}{(2n+1)(|l_1|+\cdots + |l_n| +
k_n)} \right) |I_{l_1,\ldots,l_n}|.
\label{jota-chico}
\end{equation}
The procedure of \S \ref{equiv} (using Yoccoz's equivariant family of maps) gives
subgroups of 
$\mathrm{Diff}_+^{1}([0,T_{k_1}])$ isomorphic to $H/H_n$ and generated by 
elements $A_n,B_n,C_n$ 
and $D_n$. Our next step is to estimate the $\sigma$-norm of the derivatives of
these maps.

\vspace{0.15cm}

\begin{lem} {\em If the sequence $(k_n)$ satisfy the conditions} 
\begin{equation}
\frac{(2n+1)k_n}{(2n+1)k_n - 2} \left( \frac{k_n + 1}{k_n} \right)^{2n} \leq 2, \qquad 
\left( 1-\frac{2}{(2n+1)k_n} \right) \left( \frac{k_n - 1}{k_n} \right)^{2n} \geq
\frac{1}{2},
\label{perro}
\end{equation}
\begin{equation}
2n \log(k_n) \geq \log \left( \frac{(2n+1)k_n}{(2n+1)k_n - 2} \right),
\label{tonta}
\end{equation}
and
\begin{equation}
\frac{\log(k_n)}{k_n} \left( n2^{2n+3} + \frac{32}{2n+1} \right) \leq \frac{M}{12 \pi},
\label{dos}
\end{equation}
{\em then the $\sigma$-norms of the derivatives of $A_n,B_n,C_n$ and $D_n$ 
are less than or equal to $M$ for all $n \in \mathbb{N}$.}
\end{lem}

\noindent{\bf Proof.} First of all, it is easy to verify that inequality
(\ref{jota-chico}) 
and hypothesis (\ref{perro}) imply that 
\begin{equation}
\frac{1}{2}\leq
\frac{|I_{l_1,\ldots,1+l_i,\ldots,l_n}|}{|I_{l_1,\ldots,l_i,\ldots,l_n}|}\leq 2 
\quad \mbox{ and } \quad 
\frac{1}{2}\leq
\frac{|J_{l_1,\ldots,1+l_i,\ldots,l_n}|}{|J_{l_1,\ldots,l_i,\ldots,l_n}|}\leq 2.
\label{ados}
\end{equation}
We also have
\begin{equation}
\frac{|l_1| + \cdots + |l_i| \cdots + |l_n| + k_n}{|l_1| + \cdots + |1+l_i|+\cdots
+|l_n|+k_n} 
\leq 2.
\label{gato}
\end{equation}
According to the construction of the maps and inequalities 
(\ref{ados}), the problem reduces to estimating expressions of the form 
$$\left| \frac{|I_{l_1,\ldots,1+l_i,\ldots,l_n}|}{|I_{l_1,\ldots,l_i,\ldots,l_n}|} -
1 \right| 
\frac{1}{\sigma(|I_{l_1,\ldots, l_i, \ldots,l_n}|)} 
\qquad \mbox{ and } \qquad 
\left| \frac{|J_{l_1,\ldots,1+l_i,\ldots,l_n}|}{|J_{l_1,\ldots,l_i,\ldots,l_n}|} - 1
\right| 
\frac{1}{\sigma(|J_{l_1,\ldots, l_i, \ldots,l_n}|)}.$$
Indeed, if we verify that these expressions are bounded above by $M/12 \pi$ for all
possible 
choices of sub-indices, then Lemmas \ref{simple} and \ref{pegar} will imply that the 
$\sigma$-norm of the derivatives of $A_n,B_n,C_n$ and $D_n$ 
are less than or equal to $M$.

Using the identity $s^{2n}-1 = (s-1)(s^{2n-1}+\cdots+1)$ and inequality 
(\ref{gato}) we obtain 
\begin{eqnarray*}
\left|\frac{|I_{l_1,\ldots,1+l_i,\ldots,l_n}|}{|I_{l_1,\ldots,l_i,\ldots,l_n}|} -1
\right|
&=& \left| \left( \frac{|l_1|+\cdots +|l_i|+\cdots + |l_n|+
k_n}{|l_1|+\cdots+|1+l_i|+\cdots+|l_n|+k_n} 
\right)^{2n} - 1\right|\\
&\leq&
\frac{\big| |l_i|-|1+l_i| \big|}{|l_1|+\cdots+|1+l_i|+\cdots+|l_n|+k_n} \cdot 
(2^{2n-1} + 2^{2n-2} + \cdots + 1)\\
&\leq& \frac{2^{2n}}{|l_1| + \cdots +|1+l_i|+\cdots +
|l_n|+k_n}.  
\end{eqnarray*}
Since the function $s \mapsto \log(s)/s$ is decreasing 
for $s \geq e$, by (\ref{dos}) we conclude
\begin{eqnarray*}
\left|\frac{|I_{l_1,\ldots,1+l_i,\ldots,l_n}|}{|I_{l_1,\ldots,l_i,\ldots,l_n}|} -1
\right| \frac{1}{\sigma(|I_{l_1,\ldots,l_i,\ldots,l_n}|)}
&\leq& \frac{2^{2n} \log([|l_1|+\cdots + |l_i|
 + \cdots + |l_n| + k_n]^{2n})}{|l_1|+\cdots +|1+l_i|+\cdots + |l_n|+ k_n} \\
&\leq& \frac{n2^{2n+1} \log(|l_1|+\cdots+|l_i|+\cdots
+|l_n|+k_n)}{|l_1|+\cdots+|1+l_i|+\cdots+|l_n|+k_n} \\
&\leq& \frac{n2^{2n+2}\log(k_n)}{k_n} \\
&\leq& \frac{M}{12 \pi}.
\end{eqnarray*}

Let us now deal with the case of the intervals $J_{l_1,\ldots,l_n}$. First of all, 
a straightforward computation using (\ref{jota-chico}) and (\ref{tonta}) shows that 
\begin{equation}
\sigma(|J_{l_1,\ldots,l_i,\ldots,l_n}|) \geq \sigma(|I_{l_1,\ldots,l_i,\ldots,l_n}|)
/ 2.
\label{retonta}
\end{equation}
Then, using (\ref{jota-chico}), (\ref{dos}), (\ref{retonta}), and the triangle 
inequality 
$$\left|\frac{|J_{l_1,\ldots,1+l_i,\ldots,l_n}|}{|J_{l_1,\ldots,l_i,\ldots,l_n}|} -1
\right| 
\leq \left|
\frac{|I_{l_1,\ldots,1+l_i,\ldots,l_n}|}{|I_{l_1,\ldots,l_i,\ldots,l_n}|} - 1
\right| 
+ \left| \frac{|J_{l_1,\ldots,1+l_i,\ldots,l_n}|}{|J_{l_1,\ldots,l_i,\ldots,l_n}|} - 
\frac{|I_{l_1,\ldots,1+l_i,\ldots,l_n}|}{|I_{l_1,\ldots,l_i,\ldots,l_n}|} \right|,$$
it is easy to verify that 
$$\left|\frac{|J_{l_1,\ldots,1+l_i,\ldots,l_n}|}{|J_{l_1,\ldots,l_i,\ldots,l_n}|} -1
\right| 
\frac{1}{\sigma(|J_{l_1,\ldots,l_i,\ldots,l_n}|)}$$ 
is less than or equal to
\begin{small}
$$\frac{n2^{2n+3}\log(k_n)}{k_n} + \frac{8}{(2n+1)(|l_1|+\cdots+|l_i|+\cdots+|l_n|+
k_n)} 
\cdot \frac{|I_{l_1,\ldots,1+l_i,\ldots,l_n}|}{|I_{l_1,\ldots,l_i,\ldots,l_n}|} \cdot 
\frac{1}{\sigma(|I_{l_1,\ldots,l_i,\ldots,l_n}|)},$$
\end{small}
and the value of this expression is bounded above by 
\begin{small} 
$$\frac{n2^{2n+3}\log(k_n)}{k_n} + \frac{32n}{2n+1} \cdot 
\frac{\log(|l_1|+\cdots+|l_i|+\cdots+|l_n|+k_n)}{|l_1|+\cdots+|l_i|+\cdots+|l_n|+k_n} 
\leq \frac{\log(k_n)}{k_n} \left( n2^{2n+3} + \frac{32n}{2n+1} \right) \leq
\frac{M}{12 \pi}.$$
\end{small}
This concludes the proof of the lemma.

\vspace{0.25cm}

Remark that similar computations lead to the same estimate for the 
$\sigma$-norms of $(A_n^{-1})'$, $(B_n^{-1})'$, $(C_n^{-1})'$ and $(D_n^{-1})'$. 
Thus, $A_n'$, $B_n'$, $C_n'$ and $D_n'$ converge to some $\sigma$-continuous 
functions (having $\sigma$-norm bounded above by $M$), and these functions 
are the derivatives 
of $C^{1+\sigma}$ diffeomorphisms $A,B,C$ and $D$ which generate a group 
isomorphic to $H$. Note however that this group acts on the interval $[0,T_{k_1}]$.
In order 
to obtain a group acting on $[0,1]$, we can conjugate by the affine map $g: [0,1]
\rightarrow 
[0,T_{k_1}]$, {\em i.e.} $g(x) = T_{k_1}x$. 
Since $k_1 \geq 4$ we have $T_{k_1} \leq 1$, and therefore 
this procedure does not increase $\sigma$-norms of derivatives. Indeed, for instance 
\begin{eqnarray*}
\left| \frac{(g^{-1}Ag)'(x) -
(g^{-1}Ag)'(y)}{\sigma(|x-y|)} \right| 
&=& \left| \frac{A'(g(x)) - A'(g(y))}{\sigma(|x-y|)}
\right| \\
&=& \left|
\frac{A'(g(x))-A'(g(y))}{\sigma(|g(x)-g(y)|)} \right|
\cdot \left|
\frac{\sigma(T_{k_1}|x-y|)}{\sigma(|x-y|)} \right| \\
&\leq& M.
\end{eqnarray*}
Since every (orientation preserving) 
diffeomorphism $f$ of $[0,1]$ has a point for which the derivative is equal to $1$, 
if the $\sigma$-norm of $f'$ is bounded above by $M$ then 
$$\sup_{x \in [0,1]} |f'(x) - 1| \leq M \sigma(1).$$
So, if $M$ is small, then $f$ is near the identity in the $C^{1+\sigma}$-topology. 
Finally, it is easy to construct sequences $(k_n)$ of integer positive numbers
satisfying 
(\ref{perro}), (\ref{tonta}) and (\ref{dos}). This finishes the proof of Theorem A.

\begin{rem} Note that the action we obtained is given by diffeomorphisms which are 
tangent to the identity at the end points of $[0,1]$. Therefore, gluing together 
these two end points, we obtain an action by $C^1$ diffeomorphisms of the circle. 
Using the classical procedure of suspension, this allows to construct a 
codimension-$1$ foliation (on a $3$-dimensional compact manifold) which 
is transversely of class $C^1$ and whose leaves have sub-exponential 
growth. It turns out that a countable number of leaves have polynomial 
growth and a continuum of leaves have intermediate growth, a fact that 
should be compared with \cite{hector}. In \cite{CC} one can find further 
interesting examples of codimension-$1$ foliations which are transversely 
of class $C^1$ but not $C^2$.
\end{rem}

%%%%%%%%%%%%%%%%%%%%%%%%%%%%%%%%%%%%%%%%%%%%%%%%%%%%%%%%%%%%%%%%%%%%%%%%%%%%%%%%%%%%%%%%%%%%%%%%%%%
%%%%%%%%%%%%%%%%%%%%%%%%%%%%%%%%%%%%%%%%%%%%%%%%%%%%%%%%%%%%%%%%%%%%%%%%%%%%%%%%%%%%%%%%%%%%%%%%%%%
%%%%%%%%%%%%%%%%%%%%%%%%%%%%%%%%%%%%%%%%%%%%%%%%%%%%%%%%%%%%%%%%%%%%%%%%%%%%%%%%%%%%%%%%%%%%%%%%%%%

\section{Sub-exponential growth groups of 
$\mathbf{{\it C}^{1+\alpha}}$ interval diffeomorphisms}

\hspace{0.35cm} The second part of this article 
is mainly devoted to the proof of Theorem B. 
Note that, though this theorem is stated in terms of growth, we will prove an 
{\em a priori} stronger result, namely that finitely generated subgroups of 
$\mathrm{Diff}_+^{1+\alpha}([0,1])$ {\em without free semi-groups on two generators} 
are almost nilpotent.\footnote{It seems to be unknown wether 
there exists a left orderable 
group without free semi-groups on two generators and having exponential growth.} 
For this, first recall that a classical result of J. Rosenblatt (\cite{Ros}; see 
also \cite{breuillard}) establishes that every (finitely generated) solvable 
group without free semi-groups on two generators is almost nilpotent.\footnote{It 
may be possible to give a simple proof of this result for left orderable 
groups or at least for groups of $C^1$ diffeomorphisms of the interval.} Hence, 
to prove (the {\em a priori} stronger version of) Theorem B, it suffices to 
show that finitely generated groups of $C^{1+\alpha}$ interval diffeomorphisms 
without free semi-groups on two generators are solvable.

At this point we would like to make a curious remark: for groups without free
semi-groups 
on two generators, it is no longer necessary to assume finite generation in order to 
ensure solvability. This is due to the fact that what we will show is that 
finitely generated groups of $C^{1+\alpha}$ diffeomorphisms of the interval 
without free semi-groups on two generators are solvable with solvability 
degree bounded by some constant depending only on $\alpha$.

Unfortunately, even the proof of the solvability is somehow technical. For this reason 
we decided to present the ideas in a progressive manner. In \S \ref{bp}, we recall a 
simple criterium for existence of free semi-groups arising from the theory of 
codimension-1 foliations which is closely related to the so-called 
{\em resilient leaves}. Furthermore, we discuss a useful tool also 
intimately related to all of this, namely the translation number 
homomorphism associated to an invariant Radon measure. 
In \S \ref{noH2}, we discuss the problem (already 
treated in \cite{Na-subexp}) of the embedding of 
groups with sub-exponential growth into $\mathrm{Diff}_+^2([0,1])$, 
and in \S \ref{noH}, we give a relatively short proof of the fact that, for 
any $\alpha > 0$, there is no embedding of $H$ into
$\mathrm{Diff}_+^{1+\alpha}([0,1])$ 
which is semi-conjugate to those of the first part of this work. These last 
two paragraphs can be skipped by the reader who is pressed to go into 
the proof of Theorem B. (He or she will only need Lemma \ref{cagoncito}, 
which is mostly self-contained.) However, they can serve as a useful 
guide, because it is in these two sections 
where we introduce (and explain) our main ideas: 
\S \ref{noH2} contains a general strategy of proof when it 
is possible to control distortions individually ({\em i.e.} in class $C^2$), 
whereas \S \ref{noH} contains a very simple method to get control of distortion 
in class $C^{1+\alpha}$. All these ideas are putted together in \S 
\ref{alfin}, which is divided into three subsections 
where a complete proof for Theorem B is given. Finally, 
in \S \ref{ext}, we prove the extensions of Theorem B 
for groups acting on the circle and the real line.

We would like to conclude this small introduction to the second part of 
this article by addressing the question of the uniformity of the exponential 
growth rate for groups of interval diffeomorphisms (see Remark \ref{rem-wil}):

\vspace{0.28cm}

\noindent{\bf Question.} If $\Gamma$ is a finitely generated non almost 
nilpotent subgroup of $\mathrm{Diff}^{1+\alpha}_+([0,1])$, does $\Gamma$ 
necessarily have uniform exponential growth? Is this true at least 
for non-Abelian groups of $C^2$ diffeomorphisms of the interval?

%%%%%%%%%%%%%%%%%%%%%%%%%%%%%%%%%%%%%%%%%%%%%%%%%%%%%%%%%%%%%%%%%%%%%%%%%%%%%%%%%%%%%%%%%%%%%%%%%%
%%%%%%%%%%%%%%%%%%%%%%%%%%%%%%%%%%%%%%%%%%%%%%%%%%%%%%%%%%%%%%%%%%%%%%%%%%%%%%%%%%%%%%%%%%%%%%%%%%
%%%%%%%%%%%%%%%%%%%%%%%%%%%%%%%%%%%%%%%%%%%%%%%%%%%%%%%%%%%%%%%%%%%%%%%%%%%%%%%%%%%%%%%%%%%%%%%%%%

\subsection{Crossed elements, invariant Radon measures, and translation numbers}
\label{bp}

\hspace{0.35cm} We say that two homeomorphisms of the interval $[0,1]$ are {\em
crossed} on a sub-interval $[u,v]$ if one of them fixes $u$ and $v$ and no other 
point in $[u,v]$, while the other sends $u$ or $v$ into $]u,v[$. The following 
elementary criterion for existence of free semi-groups on two 
generators is well known.

\vspace{0.1cm}

\begin{lem} {\em If a subgroup $\Gamma$ of $\mathrm{Homeo}_+([0,1])$ contains 
two crossed elements, then $\Gamma$ contains a free semi-group on two generators.} 
\label{semilibre}
\end{lem}

\noindent{\bf Proof.} Suppose that there exist $f,g$ in $\Gamma$ and an 
interval $]u,v[$ which is fixed by $f$ and contains no fixed point of $f$ 
in its interior, and such that $g(u) \! \in ]u,v[$. 
(The case where $g(v) \! \in ]u,v[$ is analogous.) Changing $f$ by its 
inverse if necessary, we may assume that $f(x) < x$ for all $x \! \in ]u,v[$. 
Let $w = g(u) \! \in ]u,v[$, and let us fix a point $z'\! \in ]w,v[$. 
Since $gf^n(u) = w$ for all $n\in\mathbb{N}$, and since $gf^n(z')$ 
converges to $w$ as $n$ tends to infinity, the map $gf^n$ has a 
fixed point in $]u,z'[$ for $n\in\mathbb{N}$ large enough. Fix such 
a $n \in \mathbb{N}$ and let $z > w$ be the infimum of the fixed points of $gf^n$ in 
$]u,v[$. For $m \in \mathbb{N}$ large enough we have $f^m(z) < w$, and so the
(positive) 
Ping-Pong Lemma applied to the restrictions of $f^m$ and $fg^n$ to $[u,v]$ shows that 
the semi-group generated by these two elements is free (see \cite{harpe}, Chapter VII).

\vspace{0.25cm}

The preceding lemma has the following important consequence: if $f$ and $g$ 
are interval homeomorphisms which generate a group without free semi-groups, 
and if $f$ has a fixed point $x_0$ which is not fixed by $g$, then the fixed 
points of $g$ immediately to the left and to the right of $x_0$ are also 
fixed by $f$. This gives a quite clear geometric picture for the action 
of a group without free semi-groups. However, it is sometimes difficult 
to use this picture without entering into rather complicated combinatorial 
discussions. In order to avoid this problem, 
there is an extremely useful tool for detecting fixed points of elements, 
namely the translation number associated to an invariant Radon measure. 
We begin by recalling a result due to J. Plante \cite{Pl3} for groups with 
sub-exponential growth, and to V. Solodov \cite{solodov} and L. Beklaryan 
\cite{Bekl} for groups without free semi-groups. We offer a short proof 
of our own for the convenience of the reader.

\vspace{0.1cm}

\begin{prop} {\em Let $\Gamma$ be a finitely generated group of homeomorphisms 
of $[0,1]$. If $\Gamma$ has no crossed elements, then $\Gamma$ preserves a 
(non trivial) Radon measure on $]0,1[$ (i.e. a measure on the Borelean 
sets which is finite on compact subsets of $]0,1[$).}
\label{radon}
\end{prop}

\noindent{\bf Proof.} If $\Gamma$ has global fixed points inside $]0,1[$ then the 
claim is obvious: the delta measure on any of such points is invariant by the 
action. Assume in what follows that the $\Gamma$-action on $]0,1[$ has no 
global fixed point, and take a finite system $\{f_1,\ldots,f_k\}$ 
of generators for $\Gamma$. We first claim that (at least) one of 
these generators does not have interior fixed points. Indeed, suppose by 
contradiction that all the maps $f_i$ have interior fixed points, and let 
$x_1 \!\in ]0,1[$ be any fixed point $f_1$. If $f_2$ 
fixes $x_1$ then letting $x_2=x_1$ we have that $x_2$ is fixed by both $f_1$ 
and $f_2$. If not, choose a fixed point $x_2 \!\in ]0,1[$ for $f_2$ 
such that $f_2$ does not fix any point between $x_1$ and $x_2$. 
Since $f_1$ and $f_2$ are non crossed on 
any interval, $x_2$ must be fixed by $f_1$. Now if $x_2$ is fixed by $f_3$ 
let $x_3=x_2$; if not, take a fixed point $x_3 \!\in ]0,1[$ for $f_3$ such that 
$f_3$ has no fixed point between $x_2$ and $x_3$. 
The same argument as before shows that $x_3$ 
is fixed by $f_1,f_2$, and $f_3$. Continuing in this way we find a 
common fixed point for all the generators $f_i$, and so a global 
fixed point for $\Gamma$, which contradicts our assumption.

Now we claim that there exists a non empty minimal invariant 
closed set for the action of $\Gamma$ inside $]0,1[$. To 
prove this consider a generator $f = f_i$ without fixed points, 
fix any point $x_0 \!\in ]0,1[$, and let $I$ be the interval 
$[x_0,f(x_0)]$ if $f(x_0) > x_0$, and $[f(x_0),x_0]$ if $f(x_0) < x_0$. In 
the family $\mathcal{F}$ of non empty closed invariant subsets of $]0,1[$ 
let us consider the order relation $\preceq$ given by \esp 
$K_1 \succeq K_2$ if \esp $K_1 \cap I \subset K_2 \cap I$. \esp 
Since $f$ has no fixed point, every orbit by $\Gamma$ must 
intersect the interval $I$, and so 
$K \cap I$ is a non empty compact set for all 
$K \in \mathcal{F}$. Therefore, we can apply Zorn Lemma to obtain a maximal 
element for the order $\preceq$, and this element is the intersection with $I$ of 
a minimal non empty closed subset of $]0,1[$ invariant by the action of $\Gamma$.

Now fix the non empty closed invariant minimal set $K$ obtained 
above. Note that its boundary $\partial K$ as well 

\vspace{-0.205cm}

\noindent as the set of its accumulation 
points $K'$ are also closed sets invariant by $\Gamma$. Because of the 
minimality of $K$, there are three possibilities:\\

\noindent -- $K' = \emptyset$: 
In this case $K$ is discrete, that is $K$ coincides 
with the set of points of a sequence $(y_n)_{n \in \mathbb{Z}}$ satisfying 
$y_{n} \!<\! y_{n+1}$ for all $n$ and without accumulation points inside $]0,1[$. 
It is then easy to see that the Radon measure \esp 
$\mu = \sum_{n \in \mathbb{Z}} \delta_{y_n}$ \esp is invariant by $\Gamma$.

\vspace{0.25cm}

\noindent -- $\partial K = \emptyset$: 
In this case $K$ coincides with the whole 
interval $]0,1[$. We claim that the action of $\Gamma$ 
is free. Indeed, if not let $[u,v]$ be an interval strictly contained in $[0,1]$ and 
for which there exists an element $g \in \Gamma$ fixing $]u,v[$ and with no fixed
point 
inside it. Since the action is minimal, there must be some $h \in \Gamma$ sending $u$ 
or $v$ inside $]u,v[$; however, this implies that $g$ and $h$ are crossed on $[u,v]$, 
contradicting our assumption. Now the action of $\Gamma$ being 
free, H\"older Theorem \cite{ghys} implies that $\Gamma$ 
is topologically conjugate to a (in this case dense) group of translations. Pulling
back 
the Lebesgue measure by this conjugacy we obtain an invariant Radon measure for the 
action of $\Gamma$.

\noindent -- $\partial K = K' = K$: 
In this case $K$ is ``locally'' a Cantor set. Collapsing to a point the closure of 
each connected component 
of the complementary set of $K$ we obtain a new (topological) open interval, and the
original 
action of $\Gamma$ induces (by semi-conjugacy) an action by homeomorphisms on this new 
interval. As in the second case, one easily checks that the induced action is free, 
and so it preserves a Radon measure. Pulling back this measure by the semi-conjugacy, 
one obtains a Radon measure on $]0,1[$ which is invariant by the original action.

\vspace{0.1cm}

\begin{rem} It is easy to see that the finite generation hypothesis is necessary 
for the preceding proposition. Note however that, during the proof, this hypothesis 
was only used to ensure the existence of a minimal non empty invariant closed 
subset of $]0,1[$, and hence it can be replaced by any other hypothesis leading 
to the same conclusion. For instance, the proposition is still true for non 
finitely generated groups containing elements without fixed points, or for 
groups having a system of generators all whose elements send each fixed point 
$x \!\in ]0,1[$ into some compact subset of $]0,1[$ (which depends on $x$).
\end{rem}

\vspace{0.02cm}

In Proposition \ref{radon} it is sometimes better to think of $\Gamma$ 
as a group of homeomorphisms of the real line. Recall that for 
(non necessarily finitely generated) groups of homeomorphisms 
of the real line preserving a Radon measure $\mu$ there is 
an associated {\em translation number} function $\tau_{\mu}: \Gamma 
\rightarrow \mathbb{R}$ defined by 

$$\tau_{\mu}(g) = \left \{ \begin{array} {l} 
\mu([x_0,g(x_0)[) \hspace{0.67cm} \mbox{ if } \esp g(x_0) > x_0,\\ 
0 \hspace{2.5cm} \mbox{ if } \esp g(x_0) = x_0,\\
- \mu([g(x_0),x_0[) \hspace{0.4cm} \mbox{ if } \esp g(x_0) < x_0, \end{array} \right.$$
where $x_0$ is any point of the real line. (One easily checks that this 
definition is independent of $x_0$.) The following properties are satisfied 
(the verification is easy; see for instance \cite{Pl3}):\\

\vspace{0.07cm}

\noindent (i) $\tau_{\mu}$ is a group homomorphism,\\

\vspace{0.07cm}

\noindent (ii) $\tau_{\mu}(g) = 0$ if and only if $g$ has fixed points; in this 
case the support of $\mu$ is contained in the set of these points,\\

\vspace{0.07cm}

\noindent (iii) $\tau_{\mu}$ is trivial if and only if $\Gamma$ has global fixed 
points ({\em i.e.} there are points which are fixed by all the elements of 
$\Gamma$).\\

\vspace{0.07cm}

Note that (iii) is a direct consequence of (ii). Moreover, (i) and 
(iii) imply the existence of global fixed points for the action of 
the derived group $\Gamma_1 = [\Gamma,\Gamma]$. 
If $\Gamma$ does not already have global fixed 
points then, according to the proof of Proposition \ref{radon}, 
there are two possibilities depending on the nature of $\mu$:\\

\vspace{0.07cm}

\noindent -- if $\mu$ has no atoms then $\Gamma$ is semi-conjugate to a group of
translations,\\

\vspace{0.07cm}

\noindent -- if $\mu$ has some atom $x_0$ then there exists an element $f \in
\Gamma$ such 
that the orbit of $x_0$ by $\Gamma$ coincides with the orbit of $x_0$ by $f$;
moreover, for 
each $g \in \Gamma$ there exists an integer $i=i(g)$ such that $g(f^n(x_0)) =
f^{n+i}(x_0)$ 
for all $n \in \mathbb{Z}$, and the homomorphism $\tau_{\mu}$ is 
given by a scalar multiple of the function $g \mapsto i(g)$.

\vspace{0.07cm}

Invariant Radon measures will be very useful because the functorial properties 
of the associated translation numbers will allow us to show that, roughly, every 
continuous action on the interval of a group  without free semi-groups on two 
generators has a level structure which is somehow similar to that of $H$.

\vspace{0.07cm}

\begin{rem} For the rest of this article one can forget about invariant Radon 
measures, and simply retain an elementary fact, 
already remarked by Salhi in 
\cite{salhi} and by Solodov in \cite{solodov}, 
and that can be checked directly: if $\Gamma$ is a (non 
necessarily finitely generated) subgroup of $\mathrm{Homeo}_+([0,1])$ 
without crossed elements, then the set of elements of $\Gamma$ having 
fixed points inside $]0,1[$ is a normal subgroup. Nevertheless, we 
think that the presentation of our arguments in terms of invariant 
Radon measures is more transparent. 
\label{rem-sold}
\end{rem}

%%%%%%%%%%%%%%%%%%%%%%%%%%%%%%%%%%%%%%%%%%%%%%%%%%%%%%%%%%%%%%%%%%%%%%%%%%%%%%%%%%%%%%%%%%%%%%%%%%
%%%%%%%%%%%%%%%%%%%%%%%%%%%%%%%%%%%%%%%%%%%%%%%%%%%%%%%%%%%%%%%%%%%%%%%%%%%%%%%%%%%%%%%%%%%%%%%%%%

\subsection{On the non embedding of the group 
$\mathbf{\it H}$ into $\mathbf{\mathrm{Diff}_+^2([0,1])}$}
\label{noH2}

\hspace{0.35cm} It is possible to give a direct and simple argument to prove that, 
for any homomorphism $\phi: H \rightarrow \mathrm{Diff}_+^2([0,1])$, the image 
$\phi(H)$ is Abelian. This proof uses only the facts that ${\tt a}^2 \in H$ 
belongs to the center of $H$ and that the elements $\tt b$, $\tt c$ and $\tt d$ 
commute between them.

Denote by $Fix_{\phi}({\tt a})$ the set of points in 
$[0,1]$ which are fixed by $\phi({\tt a}^2)$ (equivalently, 
by $\phi({\tt a})$). Fix any interval $[x_0,y_0] \subset [0,1]$ such that 
$Fix_{\phi}({\tt a}) \cap [x_0,y_0] = \{x_0,y_0\}$. Since ${\tt a}^2$ and $\tt b$ 
commute, $\phi({\tt b})^n(x_0)$ and $\phi({\tt b})^n(y_0)$ belong to
$Fix_{\phi}({\tt a})$ 
for all $n \in \mathbb{Z}$, and so they are not in $]x_0,y_0[$. If 
$\phi({\tt b})(x_0) \leq x_0$ let 
$x = \lim_{n \rightarrow \infty} \phi({\tt b})^n(x_0) \leq x_0$ 
and $y = \lim_{n \rightarrow \infty} \phi({\tt b})^{-n}(y_0) \geq y_0$. 
If $\phi({\tt b})(x_0) \geq x_0$ let 
$x = \lim_{n \rightarrow \infty} \phi({\tt b})^{-n}(x_0) \leq x_0$ 
and $y = \lim_{n \rightarrow \infty} \phi({\tt b})^{n}(y_0) \geq y_0$. Note that
both $x$ and 
$y$ belong to $Fix_{\phi}({\tt a})$. 
Kopell Lemma applied to the restrictions of $\phi({\tt a}^2)$ and 
$\phi({\tt b})$ to $[x,y]$ shows that $x=x_0$ and $y=y_0$. 
Thus the restriction of $\phi({\tt b})$ to 
$[x_0,y_0]$ is contained in the centralizer (in $\mathrm{Diff}_+^2([x_0,y_0])$) of
the restriction 
of $\phi({\tt a})^2$, 
which by Szekeres' Theorem is an Abelian group. Similar arguments can be given 
for $\phi({\tt c})$ and $\phi({\tt d})$, 
concluding that $\phi(H)$ fixes $[x_0,y_0]$ and the corresponding 
restriction is an Abelian group. On the other hand, on the set $Fix_{\phi}({\tt a})$
the 
action induced by $\phi$ factors through an action of the Abelian group generated by 
${\tt b}$, ${\tt c}$ and ${\tt d}$. This finishes the proof of the commutativity of
$\phi(H)$.

\vspace{0.2cm}

The preceding argument cannot be applied to general sub-exponential growth 
subgroups of $\mathrm{Diff}_+^2([0,1])$ because it uses commutativity 
in a way that is too strong. However, in \S \ref{bp} we saw 
that, though the center of a sub-exponential growth subgroup of 
$\mathrm{Homeo}_+([0,1])$ may be trivial, its elements have a lot of 
``dynamically commuting features''. We will see that this property still 
allows applying some of the techniques of the classical rigidity theory for 
centralizers (of $C^2$ maps). Note that this last issue already appears (as 
the main technical ingredient) in \cite{Na-solv}. Indeed, that paper contains 
the dynamical description for the actions by $C^2$ interval diffeomorphisms 
of solvable groups, for which the center is in most cases trivial...

Using the techniques introduced in \cite{Na-solv} it is proved in \cite{Na-subexp} 
that finitely generated groups of $C^2$ diffeomorphisms of the interval with 
sub-exponential growth are Abelian, thus slightly generalizing Plante-Thurston 
Theorem \cite{PT}. Actually, this result is now an easy corollary of our Theorem B. 
However, in the rest of this section we will give a complete proof of it in order 
to introduce a second criterium for proving existence of free semi-groups, 
as well as to illustrate the usefulness of invariant Radon measures and 
translation numbers to simplify many combinatorial arguments. This will also 
allow us to clarify some unclear points of the proof given in \cite{Na-subexp}. 
But before getting into the proof, let's take a look at a relevant example.

\begin{ex} The wreath product $\mathbb{Z} \wr \mathbb{Z}$ has a natural 
action on the interval. This action is obtained by identi- fying one of 
the two canonical generators of this group with a homeomorphism $f$ of $[0,1]$ 
satisfying $f(x) < x$ for all $x \! \in ]0,1[$, and the other generator with a 
homeomorphism $g$ satisfying $g(x) \neq x$ for all $x\!\in ]f(x_0),x_0[$ and $g(x)=x$ 
for all $x \in [0,1] \setminus [f(x_0),x_0]$, where $x_0$ is some point in $]0,1[$. 
(This action can be easily smoothed up to the class $C^{\infty}$; see 
\cite[Section 1]{Na-solv}.) Note that $\mathbb{Z} \wr \mathbb{Z}$ is a metabelian 
non almost nilpotent group, so according to Theorem B it must contain free 
semi-groups on two generators. Actually, the semi-group generated by $f$ 
and $g$ is free. The argument bellow is the essence of our 
second criterium for existence of free semi-groups.

\vspace{0.25cm}

\noindent{\underbar{\bf Claim:}} The semi-group generated by $f$ and $g$ is free.

\vspace{0.2cm} 

\noindent{\bf Proof.} Actually, we will prove a much more general statement. 
Suppose that $f$ and $g$ are two homeomorphisms of $[0,1]$ such that for 
some $x_0 \! \in ]0,1[$ one has $f(x_0) < x_0$, and such that $g$ fixes 
all the points of the orbit of $x_0$ by $f$. Assume moreover that for some 
interval $[u,v] \subset ]f(x_0),x_0[$ disjoint from $g(]u,v[)$ the following 
holds: for each $m \in \mathbb{N}$ the intervals $g(f^m(]u,v[))$ and 
$f^m(]u,v[)$ coincide or are disjoint, and if they coincide then $g$ 
fixes $f^n g^i([u,v])$ for every $i \in \mathbb{Z}$. In this 
situation we will show that if the set of integers
$$\mathcal{N} = \{ m \geq 0:\!\quad\! 
gf^m(]u,v[) \mbox{ and } f^m(]u,v[) \mbox{ are disjoint} \}$$
is finite, then the semi-group generated by $f$ and $g$ is free. To 
do this, let's consider two different words in positive powers of $f$ and 
$g$, and let's try to prove that they represent distinct homeomorphisms. 
After conjugacy, we may suppose that these words are of the form 
\esp \esp $W_1 = f^{n}g^{m_r}f^{n_r} \cdots g^{m_1}f^{n_1}$ \esp \esp 
and \esp \esp $W_2 = g^{q}f^{p_s}g^{q_s} \cdots f^{p_1}g^{q_1}$, 
\esp \esp where $m_j,n_j,p_j,q_j$ are positive integers, $n \geq 0$, and $q \geq 0$
(with $n > 0$ 
if $r=0$, and $q>0$ if $s=0$). Let \esp $N_1 = n_1 + \ldots + n_r + n$ \esp and 
\esp $N_2 = p_1 + \ldots + p_s$, \esp and let $m$ be the maximal value in 
$\mathcal{N}$. By the choice of $m$ one has \esp 
$W_1(f^m(u)) = f^{m + N_1}(u)$ \esp and \esp $W_2(f^m(u)) = 
f^{N_2}(g^{q_1}(f^m(u)))$. \esp However, since \esp $g^{q_1}(f^m(u)) \neq f^m(u)$ 
\esp (because $m$ belongs to $\mathcal{N}$), it is easy to verify that 
\esp $f^{N_2} (g^{q_1}(f^m(u)))$ \esp cannot be equal to 
$f^{m + N_1}(u)$. Hence \esp $W_1(f^m(u)) \neq W_2(f^m(u))$, 
\esp and so $W_1 \neq W_2$. This finishes the proof of the Claim.
\end{ex}

\vspace{0.1cm}

Now let $\Gamma$ be a finitely generated subgroup of $\mathrm{Diff}^2_+([0,1])$ 
with sub-exponential growth (or more generally, without free semi-groups on 
two generators). In order to prove that $\Gamma$ is Abelian, it is of no loss 
of generality to assume that $\Gamma$ has no global fixed point inside $]0,1[$. 
Then according to the proof of Proposition \ref{radon}, $\Gamma$ contains 
an element $f$ such that $f(x)<x$ for all $x \!\in ]0,1[$. Let $\mu$ 
be a $\Gamma$-invariant Radon measure on $]0,1[$, and let $K$ 
be the set of points inside $]0,1[$ which are fixed by all 
the elements of the first derived group $\Gamma_1 = [\Gamma,\Gamma]$. 
This set is non empty, since it contains the support of $\mu$. If $K$ 
coincides with $]0,1[$ then $\Gamma$ is Abelian. Suppose now that $K$ is 
strictly contained in $]0,1[$ and that the restriction of $\Gamma_1$ 
to each connected component of $]0,1[ \setminus K$ is free. Then by 
H\"older Theorem \cite{ghys}, the restriction of $\Gamma_1$ to every 
such connected component is Abelian, which implies that $\Gamma$ is 
metabelian. By Rosenblatt's theorem \cite{Ros}, $\Gamma$ is almost 
nilpotent, and because of Plante-Thurston Theorem, $\Gamma$ is 
almost Abelian. Finally, Szekeres' theorem implies easily that 
almost Abelian groups of $C^2$ interval diffeomorphisms are 
in fact Abelian (see \cite[Lemme 5.4]{Na-subexp}), thus 
finishing the proof in this case.

It remains the case where the action of $\Gamma_1$ on some connected 
component $I$ of the complementary set of $K$ is non free. Fix an 
element $h \in \Gamma_1$ and an interval $]u,v[$ strictly contained 
in $I$ such that $]u,v[$ is fixed by $h$ but no point inside $]u,v[$ 
is fixed by $h$. There must be some element $g \in \Gamma_1$ sending 
$]u,v[$ into a disjoint interval (contained in $I$). Indeed, if this 
is not the case then, since $\Gamma$ has no crossed elements, every 
element of $\Gamma_1$ should fix $]u,v[$, and so the points $u$ and 
$v$ would be contained in $K$, contradicting the fact that $]u,v[$ 
was strictly contained in the connected component $I$ of 
$]0,1[ \setminus K$. We will finish the proof by showing 
that the semi-group generated by $f$ and $g$ is free.

\vspace{0.35cm}
 
\noindent{\underbar{\bf{Claim:}}} There exists $N_0 \in \mathbb{N}$ 
such that $g$ fixes the interval $f^n(]u,v[)$ for every $n \geq N_0$.

\vspace{0.2cm}

\noindent{\bf Proof.} Following the proof of Kopell Lemma given in 
\cite{CClibro}, denote $I\!=]w,z[$ and fix a constant $\lambda$ such that 
$$1 < \lambda < 1 + \frac{v-u}{e^{M}(u-w)},$$ 
where $M$ is the Lipschitz constant of the function \esp $\log(f')$. \esp 
Since $g$ fixes all the points $f^n(x_0)$, its derivative at the origin must 
be equal to $1$. Let $N_0$ be such that 
\begin{equation}
g'(x) \leq \lambda \quad \mbox{and} \quad (g^{-1})'(x) \leq \lambda 
\qquad \mbox{for all } x \in f^n(I) \mbox{ and all } n \geq N_0.
\label{pichul}
\end{equation}
We will show that this $N_0$ works for the Claim (see 
Figure 3). Indeed, since $\Gamma$ has no crossed elements, if 
$g$ does not fix $f^n(]u,v[)$ then $f^n(]u,v[)$ and $g(f^n(]u,v[))$ are 
disjoint. In other words, one has \esp $g(f^n(u)) \geq f^n(v)$ \esp or \esp 
$g(f^n(v)) \leq f^n(u)$. \esp Fix $n \in \mathbb{N}$ and assume that the first 
case holds for this $n$ (for the second case just follow the same arguments 
changing $g$ by $g^{-1}$). Remark that for some 
$\bar{u} \in [u,v] \subset [w,z]$ and 
$\bar{v} \in [w,u] \subset [w,z]$ one has 
$$\frac{g(f^n(u)) - f^n(u)}{f^{n}(u) - f^{n}(w)}
\geq \frac{f^n(v) - f^n(u)}{f^n(u) - f^n(w)} = 
\frac{(f^{n})'(\bar{u})}{(f^{n})'(\bar{v})} \cdot \frac{v-u}{u-w}.$$
Since $f$ preserves $K$ and has no interior fixed points, the intervals 
in $\{f^j(I):j \in \mathbb{Z}\}$ must be pairwise disjoint. 
The well known Bounded Distortion Principle then gives 
$$\frac{(f^n)'(\bar{u})}{(f^n)'(\bar{v})} \geq \exp(-M),$$
and so 
$$\frac{g(f^n(u)) - f^n(u)}{f^{n}(u) - f^{n}(w)} 
\geq \frac{v-u}{e^{M}(u-w)} > \lambda - 1.$$ 
This implies that
$$\frac{g(f^n(u)) - f^n(w)}{f^{n}(u) - f^{n}(w)} = 1 + 
\frac{g(f^n(u)) - f^n(u)}{f^{n}(u) - f^{n}(w)} > \lambda.$$
Since $g$ fixes $f^n(w) \in K$, the left hand side member of this 
inequality is equal to $g'(x)$ for some point $x \in f^n([w,u])$. 
By (\ref{pichul}), the integer $n$ must be smaller than $N_0$, 
and this finishes the proof of the Claim.

\vspace{0.6cm}

%%%%%%%%%%%%%%%%%%%%%%%%%%%%%%%%%%%%%%%%%%%%%%%%%%%%%%%%%%%%%%%%%%%%%%%%%%%%%%%%%%%%%%%%%%%%%%
\beginpicture

\setcoordinatesystem units <1cm,1cm>

\putrule from -1 0 to 14 0

\begin{small}
\put{$w$} at 7 -0.6 
\put{$z$} at 11.5 -0.6
\put{$f^{N_0}(w)$} at 1.2 -0.6 
\put{$f^{N_0}(z)$} at 4.8 -0.6 

\putrule from 7 -1.4 to 11.5 -1.4
\putrule from 1.2 -1.4 to 4.8 -1.4
\put{$($} at 1.2 -1.4 
\put{$($} at 7 -1.4 
\put{$)$} at 4.8 -1.4 
\put{$)$} at 11.5 -1.4 
\put{$I$} at 9.25 -1.15 
\put{$f^{N_0}(I)$} at 3 -1.15  

\put{$($} at 1.2 0 
\put{$($} at 7 0 
\put{$)$} at 4.8 0 
\put{$)$} at 11.5 0
\end{small}

\begin{Large}
\put{$0$} at -1.4 -0.9
\put{$1$} at 13.6 -0.9 
\put{$($} at -1.4 0 
\put{$)$} at 13.6 0 
\end{Large}

\begin{tiny}
\put{$g(u)$} at 9.5 -0.35 
\put{$g(v)$} at 10.5 -0.35 
\put{$u$} at 7.5 -0.35 
\put{$v$} at 8.5 -0.35 
\put{$f^{N_0}(u)$} at 1.8 -0.35
\put{$f^{N_0}(v)$} at 3.2 -0.35 
\put{$($} at 7.5 0
\put{$($} at 9.5 0
\put{$($} at 1.7 0
\put{$)$} at 10.5 0
\put{$)$} at 8.5 0 
\put{$)$} at 3.3 0 
\end{tiny}

\circulararc 38 degrees from 9 1.2  
center at 6 -7.5

\circulararc 105 degrees from 10 0.25  
center at 9 -0.5  

\circulararc 180 degrees from 2.9 0.25   
center at 2.4 0.25

\put{$g$} at 9 0.4 
\put{$g$} at 2.4 0.4 
\put{$f^{N_0}$} at 6.1 1.3 
\put{$\bullet$} at 1.9 0.25 
\put{$\bullet$} at 8 0.25 
%\put{$\bullet$} at 10 0.25 
\put{$\bullet$} at 9 1.2 
%\put{$\bullet$} at 2.9 0.25 

\plot 
2.9 0.25 
2.76 0.41 / 
 
\plot 
2.9 0.25 
2.97 0.45 / 

\plot 
3 1.2 
3.32 1.2 /

\plot 
3 1.2 
3.28 1.4 /

\plot 
10 0.25 
9.8 0.35 / 

\plot 
10 0.25
9.95  0.45 / 

\put{Figure 3} at 5.67 -2

\put{} at -2 0

\endpicture

%%%%%%%%%%%%%%%%%%%%%%%%%%%%%%%%%%%%%%%%%%%%%%%%%%%%%%%%%%%%%%%%%%%%%%%%%%%%%%%%%%%%%%%%%%%%%%

\vspace{0.33cm}

Now in order to prove that the semi-group generated by 
$f$ and $g$ is free, let us consider two words of the form 
\esp \esp $W_1 = f^{n}g^{m_r}f^{n_r} \cdots g^{m_1}f^{n_1}$ \esp \esp 
and \esp \esp $W_2 = g^{q}f^{p_s}g^{q_s} \cdots f^{p_1}g^{q_1}$, 
\esp \esp where  $m_j,n_j,p_j,q_j$ are positive integers, $n \geq 0$, 
and $q \geq 0$ (with $n > 0$ if $r=0$, and $p>0$ if $s=0$). Note that 
$$\tau_{\mu}(W_1) = (n_1+\ldots+n_r+n) \esp \tau_{\mu}(f) \quad \mbox{and} 
\quad \tau_{\mu}(W_2) = (p_1+\ldots+p_s) \esp \tau_{\mu}(f).$$
So, if the values of $(n_1+\ldots+n_r+n)$ and $(p_1+\ldots+p_s)$ are 
unequal, then $W_1 \neq W_2$. Assuming that these values are equal to 
some $N \in \mathbb{N}$, the proof is finished by the following Claim.

\vspace{0.35cm}

\noindent{\underbar{\bf{Claim:}}} The elements 
$f^{-N}W_1$ and $f^{-N}W_2$ are different.

\vspace{0.15cm}

\noindent{\bf Proof.} Changing $]u,v[$ by its image by some positive iterate of $f$ 
(if necessary), we may assume that the intervals $]u,v[$ and $g(]u,v[)$ are disjoint, 
and that $g$ fixes all the intervals $f^n(]u,v[)$ for $n > 0$. Let $J$ be the 
(open) convex closure of the union \esp $\cup_{j \in \mathbb{Z}} 
\esp \esp g^j(]u,v[)$. \esp Note that $g$ fixes the interval 
$J$ and has no fixed point inside it. Now remark that
\begin{eqnarray*}
f^{-N} W_1 
&=& f^{-N} f^n g^{m_r} f^{n_r} \cdots g^{m_1} f^{n_1} \\
&=& f^{-N} f^n g^{m_r} f^{n_r} \cdots g^{m_2} f^{n_1+n_2} 
\big( f^{-n_1} g^{m_1} f^{n_1} \big) \\
&=& f^{-N} f^n g^{m_r} f^{n_r} \cdots g^{m_3} f^{n_1+n_2+n_3} 
\big( f^{-(n_1+n_2)} g^{m_2} f^{n_1+n_2} \big) 
\big( f^{-n_1} g^{m_1} f^{n_1} \big) \\
&\vdots& \\
&=& \big( f^{-(N-n)} g^{m_r} f^{N-n}\big) \cdots 
\big( f^{-(n_1+n_2)} g^{m_2} f^{n_1+n_2} \big) 
\big( f^{-n_1} g^{m_1} f^{n_1} \big)
\end{eqnarray*}
and
\begin{eqnarray*}
f^{-N} W_2 
&=& f^{-N} g^{q} f^{p_s} g^{q_s} \cdots f^{p_1} g^{q_1} \\
&=& f^{-N} g^{q} f^{p_s} g^{q_s} \cdots f^{p_3} g^{q_3} f^{p_1+p_2} 
\big( f^{-p_1} g^{q_2} f^{p_1} \big) g^{q_1} \\
&\vdots& \\
&=& \big(f^{-N} g^q f^N \big) \cdots 
\big( f^{-p_1} g^{q_2} f^{p_1} \big) g^{q_1}.
\end{eqnarray*}
Since $\Gamma$ has no crossed elements, and since all the maps 
$$\big( f^{-(N-n)} g^{m_r} f^{N-n}\big), \esp \ldots, \esp  
\big( f^{-(n_1+n_2)} g^{m_2} f^{n_1+n_2} \big), \esp 
\big( f^{-n_1} g^{m_1} f^{n_1} \big) \quad \mbox{and}
\quad \big(f^{-N} g^q f^N \big), \esp \ldots, \esp 
\big( f^{-p_1} g^{q_2} f^{p_1} \big)$$
have fixed points inside $J$, they must fix the interval $J$. 
On the other hand, $g^{q_1}$ fixes $J$ but has no fixed point 
inside it. Therefore, if $\nu$ is any Radon measure on $J$ which 
is invariant by the group generated by (the restrictions to 
$J$ of) all those maps (including $g^{q_1}$), then 
\esp $\tau_{\nu}(f^{-N}W_1) = 0$ \esp and \esp 
$\tau_{\nu}(f^{-N}W_2) = \tau_{\nu}(g^{q_1}) \neq 0$, 
\esp and this shows that $f^{-N} W_1 \neq f^{-N} W_2$.

\vspace{0.05cm}

\begin{rem} It is not difficult to adapt the arguments of this 
section to prove a similar result for groups of diffeomorphisms of 
the interval $[0,1]$ having derivatives with finite total variation.
\label{clog}
\end{rem}

%%%%%%%%%%%%%%%%%%%%%%%%%%%%%%%%%%%%%%%%%%%%%%%%%%%%%%%%%%%%%%%%%%%%%%%%%%%%%%%%%%%%%%%%%%%%%%%
%%%%%%%%%%%%%%%%%%%%%%%%%%%%%%%%%%%%%%%%%%%%%%%%%%%%%%%%%%%%%%%%%%%%%%%%%%%%%%%%%%%%%%%%%%%%%%%

\subsection{On the non embedding of the group 
$\mathbf{\it H}$ into $\mathrm{Diff}_+^{1+\alpha}([0,1])$}
\label{noH}

\hspace{0.35cm} Like in the beginning of \S \ref{noH2}, one 
can also give a direct proof of the fact that the action of $H$ 
of \S \ref{cantor} is not the semi-conjugate of an action by $C^{1+\alpha}$ 
diffeomorphisms for any $\alpha > 0$. For this, it suffices to note that 
inside $H$ there are ``a lot'' of commuting elements, and then apply 
the Generalized Kopell Lemma ({\em i.e.} Th\'eor\`eme B of \cite{DKN}). 
But again, we would like to give a proof of this fact which does 
not use commutativity in an essential way. Since we are prescribing 
the combinatorial structure for the dynamics of the action (up 
to topological semi-conjugacy), the main difficulty for this 
will consist in getting good control of distortion 
estimates in class $C^{1+\alpha}$. What follows 
is much inspired by \S 1.2 of \cite{DKN}; for instance, 
the lemma below appears as Lemme 1.3 therein. For the reader's 
convenience, and because of its simplicity, we recall its proof.

\vspace{0.16cm}

\begin{lem} {\em Let $h$ be any $C^{1+\alpha}$ diffeomorphism of a closed interval
$[u,v]$. 
If $M$ denotes the $\alpha$-H\"older constant for $h'$, then for every $x \in [u,v]$ 
one has \esp \esp $|h(x) - x| \leq M \esp |v-u|^{1+\alpha}$.}
\label{cagoncito}
\end{lem}

%\vspace{-0.6cm}

\noindent{\bf Proof.} By the Mean Value Theorem, there exist 
some points \esp $y \in [u,x]$ and $z \in [u,v]$ such that
$$h'(y) = \frac{h(x) - u}{x - u} \qquad \mbox{ and } 
\qquad h'(z) = \frac{h(v)-h(u)}{v-u} = 1.$$
Since \esp \esp 
$|h'(y) - 1| = |h'(y) - h'(z)| \leq M \esp |y - z|^{\alpha} 
\leq M \esp |v-u|^{\alpha}$, \esp \esp this gives 
$$|h(x) - x| = |x - u| \esp |h'(y) - 1| \leq M |v-u|^{1+\alpha}.$$

\vspace{0.25cm}

Now let us suppose that for some $\alpha > 0$ there exists an embedding 
$H \hookrightarrow \mathrm{Diff}_+^{1+\alpha}([0,1])$ which is semi-conjugate 
to that of Example \ref{emb-lip}. Fix the smallest positive integer $k$ such that 
\begin{equation}
\alpha \esp (1 + \alpha)^{k-2} \geq 1.
\label{defk}
\end{equation}
Note that viewing $k=k(\alpha)$ as a function of $\alpha$ we have 
that $k(\alpha) \rightarrow \infty$ as $\alpha \rightarrow 0$. For each 
$(l_1,\ldots,l_k) \in \mathbb{Z}^k$ denote by $L_{l_1,\ldots,l_k}$ 
the preimage under the semi-conjugacy of the corresponding interval
$I_{l_1,\ldots,l_k}$. 
We leave to the reader the easy task of verifying the existence of elements
$h_1,\ldots,h_k$ 
in $H$ such that, for each $j \in \{1,\ldots,k\}$ and each $(l_1,\ldots,l_k) \in
\mathbb{Z}^k$, 
\begin{equation}
h_j (L_{l_1,\ldots,l_{j-1},l_j,l_{j+1},\ldots,l_k}) 
= L_{l_1,\ldots,l_{j-1},l_j-2,l_{j+1},\ldots,l_k}.
\label{call}
\end{equation}
For instance, one can take $h_1 = {\tt a}^{-2}$, $h_2 = {\tt b}^{-2}$, 
$h_3 = {\tt a}^{-1} {\tt b}^{-2} {\tt a}$, etc. The contradiction 
is then given by the following general proposition.

\vspace{0.18cm}

\begin{prop} {\em Given an integer $k \geq 3$ let 
$\{ L_{l_1,\ldots,l_k}: (l_1,\ldots,l_k) \in \mathbb{Z}^k \}$ be a family 
of closed intervals with disjoint interiors and disposed inside $[0,1]$ respecting 
the (direct) lexicographic order, that is, $L_{l_1,\ldots,l_k}$ is to 
the left of $L_{l_1',\ldots,l_k'}$ if and only if $(l_1,\ldots,l_k)$ is 
lexicographically smaller than 
$(l_1',\ldots,l_k')$. Let $h_1,\ldots,h_k$ be $C^1$ diffeomorphisms of 
$[0,1]$ such that for each $j \in \{1,\ldots,k\}$ and each 
$(l_1,\ldots,l_k) \in \mathbb{Z}^k$ one has} 
\begin{equation}
h_j (L_{l_1,\ldots,l_{j-1},l_{j},\ldots,l_k}) = 
L_{l_1,\ldots,l_{j-1},l_{j}',\ldots,l_k'} \quad 
\mbox{for some} \quad (l_j',l_{j+1}',\ldots,l_{k}') \in \mathbb{Z}^{k-j+1} 
\mbox{ satisfying} \quad l_j' \neq l_j.
\label{compli}
\end{equation}
{\em If $\alpha > 0$ is such that $k \geq k(\alpha)$, then $h_1,\ldots,h_{k-1}$ 
cannot be simultaneously contained in $\mathrm{Diff}_+^{1+\alpha}([0,1])$.}
\label{cagon}
\end{prop}

\noindent{\bf Proof.} Suppose by contradiction that $h_1,\ldots,h_{k-1}$ 
are $C^{1+\alpha}$ diffeomorphisms of $[0,1]$. Let $M \geq 1$ 
be a simultaneous $\alpha$-H\"older constant for $h_2,\ldots,h_{k-1}$ 
as well as for $\log(h_1')$, and let $\bar{M}$ be defined by
$$\log(\bar{M}) = 
\big( 1 + \alpha [1 + (1+\alpha) + 
(1+\alpha)^2 + \ldots + (1+\alpha)^{k-3}] \big) \esp \log(M).$$
Denote $[a_k,b_k] = L_{0,\ldots,0}$, 
and for $i \in \{k-1,\ldots,1\}$ define by induction the interval 
$[a_i,b_i]$ as being the closed convex closure of the union 
$$\bigcup_{n \in \mathbb{Z}} \esp h_{i+1}^n \big( [a_{i+1},b_{i+1}] \big).$$
Let $\lambda$ be a constant such that 
$$1 < \lambda < 1 + \frac{b_k - a_k}{e^{\bar{M}} \esp (b_k - a_{k-1})}.$$
The derivative of $h_k$ must be equal to $1$ at the accumulation points 
of the intervals $L_{i_1,0,\ldots,0}$ as $|i_1| \rightarrow \infty$. Therefore, 
we can fix a (large) positive integer $n$ such that $h_k' (x) \leq \lambda$ 
and $(h_k^{-1})'(x) \leq \lambda$ for all $x \in h_1^n ([a_{k-1},b_{k-1}])$. The 
interval $h_1^n ([a_k,b_k])$ is of the form $L_{l_1',\ldots,l_k'}$ for 
some $(l_1',\ldots,l_k') \in \mathbb{Z}^k$ satisfying $l_1' = l_1 + nd$ 
for some fixed $d \neq 0$. Hence, by the hypothesis (\ref{compli}), the 
intervals $h_1^n(]a_k,b_k[)$ and $h_k (h_1^n(]a_k,b_k[))$ are disjoint. 
Assume for instance that $h_k (h_1^n(a_k)) \geq h_1^n(b_k)$ (for the case 
where $h_k (h_1^n(b_k)) \leq h_1^n(a_k)$ just follow the same arguments 
changing $h_k$ by $h_k^{-1}$). Then for some $u \in [a_k,b_k] \subset 
[a_{k-1},b_{k-1}]$ and $v \in [a_{k-1},a_k] \subset [a_{k-1},b_{k-1}]$ one has 
\begin{equation}
\frac{h_k(h_1^{n}(a_k)) - h_1^{n}(a_k)}{h_1^{n}(a_k) - h_1^{n}(a_{k-1})} \geq 
\frac{h_1^{n}(b_k) - h_1^{n}(a_k)}{h_1^{n}(a_k) - h_1^{n}(a_{k-1})} = 
\frac{(h_1^{n})'(u)}{(h_1^{n})'(v)} \cdot \frac{b_k-a_k}{a_k - a_{k-1}}.
\label{mani}
\end{equation}
Note that 
\begin{eqnarray*}
\left| \log \Big( \frac{(h_1^{n})'(u)}{(h_1^{n})'(v)} \Big) \right| 
&=& \left| \log \Big( \prod_{j=0}^{n - 1} \frac{h_1'(h_1^j(u))}{h_1'(h_1^j(v))}
\Big) \right| \\
&\leq& \sum_{j=0}^{n - 1} \big| \log(h_1'(h_1^j(u))) - \log(h_1'(h_1^j(v))) \big| \\
&\leq& M \esp \sum_{j=0}^{n - 1} \big| h_1^j(u) - h_1^j(v) \big|^{\alpha} \\
&\leq& M \esp \sum_{j=0}^{n - 1} \big| h_1^j(b_{k-1}) - h_1^j (a_{k-1})
\big|^{\alpha}. 
\end{eqnarray*}
Now using the combinatorial hypothesis (\ref{compli}) and applying 
Lemma \ref{cagoncito} to \esp $h_{k-1}, \ldots, h_{2}$ \esp one obtains 

\begin{eqnarray*} 
\big| h_1^j(b_{k-1}) - h_1^j(a_{k-1}) \big| 
&\leq& M \esp \big| h_1^j(b_{k-2}) - h_1^j(a_{k-2}) \big|^{1+\alpha} \\
&\leq& M \esp \big( M \esp \big| h_1^j(b_{k-3}) - h_1^j(a_{k-3}) \big|^{1+\alpha}
\big)^{1+\alpha} 
= M^{1 + (1+\alpha)} \esp \big| h_1^j(b_{k-3}) - h_1^j(a_{k-3}) \big|^{(1+\alpha)^2} \\
&\vdots& \\
&\leq& M^{1 + (1+\alpha) + (1+\alpha)^2 + \ldots + (1+\alpha)^{k-3}} 
\esp \big| h_1^j(b_{1}) - h_1^j(a_{1}) \big|^{(1+\alpha)^{k-2}}.
\end{eqnarray*}
We then deduce 
\begin{equation}
\left| \log \Big( \frac{(h_1^{n})'(u)}{(h_1^{n})'(v)} \Big) \right| 
\leq \bar{M} \esp \sum_{j=0}^{n-1} 
\big| h_1^j(b_{1}) - h_1^j(a_{1}) \big|^{\alpha(1+\alpha)^{k-2}}.
\label{chita}
\end{equation}
The hypothesis $k \geq k(\alpha)$ implies that $\alpha (1+\alpha)^{k-2} \geq 1$. 
Therefore, since the intervals $h_1^j (]a_{1},b_{1}[)$ are two by two 
disjoint, the right hand side expression of (\ref{chita}) is bounded 
by $\bar{M}$, which implies that 
$$\frac{(h_1^{n})'(u)}{(h_1^{n})'(v)} \geq \frac{1}{e^{\bar{M}}}.$$
Introducing this inequality into (\ref{mani}) we get 
$$\frac{h_k(h_1^{n}(a_k)) - h_1^{n}(a_k)}{h_1^{n}(a_k) - h_1^{n}(a_{k-1})} 
\geq \frac{b_k - a_k}{e^{\bar{M}} \esp (a_k - a_{k-1})},$$
and summing 1 to both members this gives
$$\frac{h_k(h_1^{n}(a_k)) - h_1^{n}(a_{k-1})}{h_1^{n}(a_k) - h_1^{n}(a_{k-1})} 
\geq 1 + \frac{b_k - a_k}{e^{\bar{M}} \esp (a_k - a_{k-1})} > \lambda.$$
But since $h_k(h_1^{n}(a_{k-1})) = h_1^{n}(a_{k-1})$, 
the left hand side member in this inequality 
is equal to $h_k'(x)$ for some point $x \in h_1^{n}([a_{k-1},b_{k-1}])$, 
and so by our choice of $n$ it is less than or equal to $\lambda$. 
This contradiction finishes the proof.

\begin{rem} Using the methods of \cite{DKN}, it is possible to prove that 
the preceding proposition is still true for $k > 1 + 1/\alpha$; moreover, 
this regularity is sharp, in the sense that for every $\alpha$ such that 
$k < 1 + 1/\alpha$ there exist $C^{1+\alpha}$ counter-examples.
\label{otimo}
\end{rem}

%%%%%%%%%%%%%%%%%%%%%%%%%%%%%%%%%%%%%%%%%%%%%%%%%%%%%%%%%%%%%%%%%%%%%%%%%%%%%%%%%%%%%%%%%%%%%%
%%%%%%%%%%%%%%%%%%%%%%%%%%%%%%%%%%%%%%%%%%%%%%%%%%%%%%%%%%%%%%%%%%%%%%%%%%%%%%%%%%%%%%%%%%%%%%

\subsection{Proof of Theorem B}
\label{alfin}

\hspace{0.35cm} In what follows $\Gamma$ will be supposed to be a finitely generated 
subgroup of $\mathrm{Diff}_+^{1+\alpha}([0,1])$ without free semi-groups on two 
generators. As explained at the beginning of the second part of this 
article, to prove that $\Gamma$ is almost nilpotent it suffices to 
show that $\Gamma$ is solvable; in fact, we will prove that the 
corresponding degree of solvability is bounded by $k=1+k(\alpha)$, where 
$k(\alpha)$ is the smallest integer satisfying the inequality 
(\ref{defk}).\footnote{According to \cite{DKN}, it is quite possible 
that the solvability degree is bounded by any integer $k$ 
satisfying $k > 1 + 1/\alpha$ (compare with Remark \ref{otimo}). 
An analogous remark applies to \S \ref{ext}.} Note that, because of this 
uniform bound in terms of $\alpha$, we may (and we will) assume that the action 
of $\Gamma$ on $]0,1[$ has no global fixed point. By the first part 
of the proof of Proposition \ref{radon}, every finite system of generators of 
$\Gamma$ contains an element without interior fixed points. Let us fix once and for all 
an element $h_1 \in \Gamma$ such that $h_1(x) < x$ for all $x \! \in ]0,1[$.

%%%%%%%%%%%%%%%%%%%%%%%%%%%%%%%%%%%%%%%%%%%%%%%%%%%%%%%%%%%%%%%%%%%%%%%%%%%%%%%%%%%%%%%%%%%%%%
%%%%%%%%%%%%%%%%%%%%%%%%%%%%%%%%%%%%%%%%%%%%%%%%%%%%%%%%%%%%%%%%%%%%%%%%%%%%%%%%%%%%%%%%%%%%%%

\subsubsection{Capturing minimal levels}
\label{changing}

\hspace{0.35cm} Let us denote by $\Gamma_{\star}$ the subgroup of $\Gamma$
consisting of 
the elements having fixed points inside $]0,1[$ . Note that the fact that
$\Gamma_{\star}$ 
is a subgroup follows from that $\Gamma_{\star}$ can be identified with the kernel
of the 
translation number homomorphism $\tau_{\mu}$ associated to any fixed
$\Gamma$-invariant 
Radon measure $\mu$ on $]0,1[$. (This can also be directly checked, and holds more 
generally for groups without crossed elements which are non necessarily finitely 
generated; see Remark \ref{rem-sold}). In particular, it contains the first 
derived group $\Gamma_1$. For each $h \in \Gamma_{\star}$ distinct from 
the identity denote by $\mathcal{I}(h)$ the family of open intervals 
$I$ such that $h$ fixes $I$ but no point inside $I$ is fixed by $h$. 
For $I \in \mathcal{I}(h)$ put
$$\mathcal{N}(h,I) = \{n \geq 0: \esp h_1^{-n} h h_1^n \mbox{ has no fixed point 
inside } I \}.$$
Note that $n\!=\!0$ belongs to $\mathcal{N}(h,I)$. 
Moreover, since $\Gamma$ does not have crossed elements, if $n$ is contained 
in $\mathcal{N}(h,I)$ then the intervals $h_1^{-n} h h_1^n (I)$ and $I$ 
either coincide or are disjoint. For each $n \in \mathcal{N}(h,I)$ 
define $I_n$  as being the (open) convex closure of the union \esp 
$\cup_{j \in \mathbb{Z}} \esp h_1^{-n} h^j h_1^n (I)$. \esp Now let 
us consider the preorder\footnote{Recall that a {\em preoder} is a 
relation which is reflexive and transitive, but not necessarily 
antisymmetric.} relation $\preceq$ on $\mathcal{N}(h,I)$ defined 
by $m \preceq n$ if $I_m \subset I_n$. Again the non existence of crossed 
elements implies that for $m \prec n$ the geometric picture is as in Figure 
4. Our first task is to prove that there are integers $n$ for which $I_n$ 
is maximal. We would like to point out that there is a somehow related 
issue in the classical level theory for codimension-1 foliations which concerns 
the existence of local minimal sets (see \cite[Theorem 8.1.8]{CClibro}). 
However, in that context this property is established using $C^2$ control 
of distortion estimates which are no longer available in the $C^{1+\alpha}$ 
case. This is the main reason why our argument is so different: we will 
show that ``there must exist some minimal level" because of the absence 
of free semi-groups on two generators inside $\Gamma$. (This should be 
compared with Lemma 2.3 of \cite{PT}.)

\vspace{-0.25cm}

%%%%%%%%%%%%%%%%%%%%%%%%%%%%%%%%%%%%%%%%%%%%%%%%%%%%%%%%%%%%%%%%%%%%%%%%%%%%%%%%%%%%%%%%%%%%%%
\beginpicture
\setcoordinatesystem units <1cm,1cm>

\putrule from -1.5 0 to 14 0

\begin{small}
\put{$h_1^{-m}hh_1^m$} at 6.2 1.08 
\put{$($} at 4 0 
\put{$)$} at 8.4 0
\put{$I_m$} at 7.9 -0.5
\end{small}

\begin{tiny}
\put{$($} at 4.5 0
\put{$($} at 5.7 0
\put{$($} at 6.9 0
\put{$)$} at 5.1 0 
\put{$)$} at 6.3 0 
\put{$)$} at 7.5 0
\end{tiny}

\begin{Large}
\put{$($} at -2.1 0 
\put{$)$} at 13.4 0 
\put{$I_n$} at 12.6 0.4 
\end{Large}

\put{$\cdots$} at 3.8 0.25 
\put{$\cdots$} at 7.4 0.25
\put{$\cdots$} at -1.68 -0.45 
\put{$\cdots$} at 12.2 -0.45
\put{$\cdots$} at -1.21 -0.45 
\put{$\cdots$} at 12.67 -0.45
\put{$h_1^{-n}hh_1^n$} at 1.8 -1.5 
\put{$h_1^{-n}hh_1^n$} at 10.5 -1.3 
\putrule from 3.4 -0.7 to 7.8 -0.7
\put{$($} at 3.4 -0.7 
\put{$)$} at 7.8 -0.7
\put{$($} at 8.7 0
\put{$($} at -0.5 0
\put{$)$} at 2.3 0 
\put{$)$} at 11.5 0
\put{$\bullet$} at 1 -0.5
\put{$\bullet$} at 6.5 -1.2 

\circulararc 38 degrees from 1 -0.5   
center at 4 5 

\circulararc -38 degrees from 10.5 -0.5   
center at 7.5 5.

\circulararc 180 degrees from 5.4 0.25 
center at 4.9 0.25

\circulararc 180 degrees from 6.8 0.25 
center at 6.3 0.25

\put{$\bullet$} at 4.4 0.25 
\put{$\bullet$} at 5.8 0.25

\plot 
5.4 0.25 
5.5 0.46 / 
 
\plot 
5.4 0.25 
5.22 0.43 / 

\plot 
6.8 0.25 
6.9 0.46 / 
 
\plot 
6.8 0.25 
6.62 0.43 / 

\plot 
5.02 -1.18
4.9 -1.31 /

\plot 
5.02 -1.18
4.85 -1.12 /

\plot 
10.5 -0.5  
10.29 -0.49 / 

\plot 
10.5 -0.5  
10.38 -0.68 / 

\put{Figure 4} at 5.9 -2

\put{} at -2.3 2.2

\endpicture

%%%%%%%%%%%%%%%%%%%%%%%%%%%%%%%%%%%%%%%%%%%%%%%%%%%%%%%%%%%%%%%%%%%%%%%%%%%%%%%%%%%%%%%%%%%%%%

\vspace{0.23cm}

\begin{lem} {\em For every non trivial element $h$ in $\Gamma_{\star}$ 
and every interval $I \in \mathcal{I}(h)$, the preorder relation 
$\preceq$ on $\mathcal{N}(h,I)$ has a maximal element.}
\label{max}
\end{lem}

\noindent{\bf Proof.} Suppose by contradiction that there is no maximal element for 
$\preceq$, and fix an increasing sequence $(k_i)$ of non negative integers such that 
\esp $k_{i} \succ n$ \esp holds for every $n < k_i$. We will use this sequence to 
prove that the semi-group generated by $h_1^{-1}$ and $h$ is free, thus providing a 
contradiction. So let us consider two different words in positive powers of $h_1^{-1}$ 
and $h$, and let's try to prove that they do not represent the same element of 
$\Gamma$. After conjugacy, we may assume that these words are of the form \esp 
\esp $W_1 = h_1^{-n} h^{m_r} h_1^{-n_r} \cdots h^{m_1} h_1^{-n_1}$ \esp \esp 
and \esp \esp $W_2 = h^{q} h_1^{-p_s} h^{q_s} \cdots h_1^{-p_1} h^{q_1}$, 
\esp \esp where $m_j,n_j,p_j,q_j$ are positive integers, $n \geq 0$, and 
$q \geq 0$ (with $n > 0$ when $r=0$, and $q > 0$ when $s=0$).

We first claim that if the integers \esp $N_1 = n_1 + \ldots + n_r + n$ 
\esp and \esp $N_2 = p_1 + \ldots + p_s$ \esp are different, then $W_1$ 
and $W_2$ are distinct elements. Indeed, since $h$ is in $\Gamma_{\star}$ 
one has $\tau_{\mu} (h) = 0$, and so \esp $\tau_{\mu}(W_1) = -N_1 \esp
\tau_{\mu}(h_1)$ 
\esp and \esp $\tau_{\mu}(W_2) = -N_2 \esp \tau_{\mu}(h_1)$. On the other hand, 
since $h_1$ has no fixed point inside $]0,1[$ one has $\tau_{\mu}(h_1) \neq 0$. 
Therefore, \esp $\tau_{\mu}(W_1) \neq \tau_{\mu}(W_2)$ \esp 
when \esp $N_1 \neq N_2$.

Assume in the rest of the proof that $N_1$ and $N_2$ are equal, and denote by $N$ 
their common value. Fix an integer $i \in \mathbb{N}$ such that $k_i \geq N$. We 
will prove that $h_1^{N} W_1$ and $h_1^{N} W_2$ are different. For this, note that
\begin{eqnarray*}
h_1^{N} W_1 
&=& h_1^{N} h_1^{-n} h^{m_r} h_1^{-n_r} \cdots h^{m_1} h_1^{-n_1} \\
&=& h_1^{N} h_1^{-n} h^{m_r} h_1^{-n_r} \cdots h^{m_2} h_1^{-(n_1+n_2)} 
\big( h_1^{n_1} h^{m_1} h_1^{-n_1} \big) \\
&=& h_1^{N} h_1^{-n} h^{m_r} h_1^{-n_r} \cdots h^{m_3} h_1^{-(n_1+n_2+n_3)} 
\big( h_1^{n_1+n_2} h^{m_2} h_1^{-(n_1+n_2)} \big) \esp 
\big( h_1^{n_1} h^{m_1} h_1^{-n_1} \big) \\
&\vdots& \\
&=& \big( h_1^{N-n} h^{m_r} h_1^{-(N-n)} \big) \cdots  
\big( h_1^{n_1+n_2} h^{m_2} h_1^{-(n_1+n_2)} \big) \esp 
\big( h_1^{n_1} h^{m_1} h_1^{-n_1} \big),
\end{eqnarray*} 
and 
\begin{eqnarray*}
h_1^{N} W_2 
&=& h_1^{N} h^{q} h_1^{-p_s} h^{q_s} \cdots h_1^{-p_2}h^{q_2}h_1^{-p_1}h^{q_1} \\
&=& h_1^{N} h^{q} h_1^{-p_s} h^{q_s} \cdots 
h^{q_3} h_1^{-(p_1+p_2)} \big( h_1^{p_1} h^{q_2} h_1^{-p_1} \big) h^{q_1} \\
&\vdots& \\
&=& \big( h_1^{N} h^q h_1^{-N} \big) \esp 
\big( h_1^{N - p_s} h^{q_s} h_1^{-(N - p_s)} \big) \cdots 
\big( h_1^{p_1} h^{q_2} h_1^{-p_1} \big) h^{q_1}.
\end{eqnarray*}
Now from the facts that $k_i \geq N$ and that $k_i \succ n$ 
for all $n<k_i$ it follows easily that all the maps 
$$\big( h_1^{N-n} h^{m_r} h_1^{-(N-n)} \big), \ldots,   
\big( h_1^{n_1+n_2} h^{m_2} h_1^{-(n_1+n_2)} \big), \esp 
\big( h_1^{n_1} h^{m_1} h_1^{-n_1} \big),$$
as well as 
$$\quad \big( h_1^{N} h^q h_1^{-N} \big), \esp 
\big( h_1^{N - p_s} h^{q_s} h_1^{-(N - p_s)} \big), \ldots, 
\big( h_1^{p_1} h^{q_2} h_1^{-p_1} \big), \esp h^{q_1},$$
fix the interval $h_1^{k_i}(I_{k_i})$, but only the last one (namely 
$h^{q_1}$) has no fixed point inside it. This implies that if $\nu$ 
is a Radon measure on $h_1^{k_i}(I_{k_i})$ which is invariant by the group 
generated by all of them, then
$$\tau_{\nu} \big( h_1^{N} W_1 \big) = 0 \qquad \mbox{and} \qquad 
\tau_{\nu} \big( h_1^{N} W_2 \big) = \esp \tau_{\nu}(h^{q_1}) \neq 0.$$
In particular, $h_1^NW_1 \neq h_1^NW_2$, and so $W_1 \neq W_2$, 
which finishes the proof of the lemma.

\vspace{0.35cm}

Now let us define the {\em star operation} of ``minimization of levels'' 
as follows: for each non trivial $h \in \Gamma_{\star}$ and each 
$I \in \mathcal{I}(h)$, fix an integer $n \in \mathcal{N}(h,I)$ 
which is maximal for the preorder relation $\preceq$, and denote 
$$h^* = h^*_I = h_1^{-n} h h_1^n \qquad \mbox{and} \qquad I^*(h,I) = I_n.$$
Note that the choice of $n$ (and hence the choice of $h^*$) 
is not necessarily unique, but the interval $I^*(h,I)$ does 
not depend on it. For simplicity, we will assume in addition 
that $h^* = h$ when $n=0$ belongs to $\mathcal{N}(h,I)$. 
Remark finally that $h^*$ is always in $\Gamma_{\star}$, 
and $I^*(h,I)$ is an element of $\mathcal{I}(h^*)$. 
Moreover, it follows easily from the definitions 
that for all $m \geq 0$ the element $h$ fixes 
the interval $h_1^m(I^*(h,I))$.

%%%%%%%%%%%%%%%%%%%%%%%%%%%%%%%%%%%%%%%%%%%%%%%%%%%%%%%%%%%%%%%%%%%%%%%%%%%%%%%%%%%%%%%%%%%%%%
%%%%%%%%%%%%%%%%%%%%%%%%%%%%%%%%%%%%%%%%%%%%%%%%%%%%%%%%%%%%%%%%%%%%%%%%%%%%%%%%%%%%%%%%%%%%%%

\subsubsection{On the structure of minimal levels}
\label{diff}

\hspace{0.35cm} The commutator of a finitely generated group is not necessarily 
finitely generated, even in the case of groups of intermediate growth (it 
is however the case for nilpotent groups). This is a small algebraic 
difficulty for the proof of the following lemma. 

\vspace{0.1cm} 

\begin{lem}{\em If $\Gamma$ is not solvable with degree of solvability less than 
or equal to $k$, then there exist elements $h_2,\ldots,h_k$ in $\Gamma_{\star}$, 
and open intervals $I_1,\ldots,I_{k-1}$, such that $I_j$ is strictly contained 
in $I_i$ for $i \!<\! j$, the interval $I_{i-1}$ coincides with $I^*(h_{i},I_{i})$ 
for each $i \in \{2,\ldots,k\}$, and one has $h_i = (h_i)^*_{I_i}$ for each 
$i \in \{2,\ldots,k-1\}$.}
\label{diferencia}
\end{lem}

\vspace{0.1cm}

\noindent{\bf Proof.} Let us start by considering a non trivial element $f_k$ 
in the $k$-derived group $\Gamma_k$ of $\Gamma$, and let us fix an interval $J$ 
in $\mathcal{I}(f_k)$. We then let $h_k = f_k^*$ and $I_{k-1} = I^*(f_k,J)$. 
We claim that there exists an element $f_{k-1} \in \Gamma_{k-1}$ such that 
$I_{k-1}$ and $f_{k-1}(I_{k-1})$ are disjoint. Indeed, if this is not 
true then by Lemma \ref{semilibre} one has \esp $h(I_{k-1}) = I_{k-1}$ 
\esp for every $h \in \Gamma_{k-1}$. Now $h_k$ can be written as a product 
$[g_1,g_2] \cdots [g_{2n-1},g_{2n}]$, where the elements $g_1,\ldots,g_{2n}$ 
belong to $\Gamma_{k-1}$. The group $\Gamma_{g_1,\ldots,g_{2n}}$ generated by 
them fixes the interval $I_{k-1}$, and since it contains no free semi-group 
on two generators, it preserves a Radon measure on $I_{k-1}$. By 
properties (i) and (iii) of \S \ref{bp}, the derived group 
$[\Gamma_{g_1,\ldots,g_{2n}},\Gamma_{g_1,\ldots,g_{2n}}]$ 
has global fixed points in $I_{k-1}$. In particular, the 
element $h_k \in [\Gamma_{g_1,\ldots,g_{2n}},\Gamma_{g_1,\ldots,g_{2n}}]$ 
has fixed points inside $I_{k-1}$, which contradicts the fact that $I_{k-1}$ 
belongs to $\mathcal{I}(h_k)$.

With respect to $I_{k-1} \in \mathcal{I}(f_{k-1})$ we can 
define $h_{k-1} = f_{k-1}^*$ and $I_{k-2} = I^*(f_{k-1},I_{k-1})$. 
As above, one can prove the existence of an element $f_{k-2} \in \Gamma_{k-2}$ 
such that $I_{k-2}$ and $f_{k-2}(I_{k-2})$ are disjoint, and then with respect 
to $I_{k-2} \in \mathcal{I}(f_{k-2})$ we let $h_{k-2} 
= f_{k-2}^*$ and $I_{k-3} = I^*(f_{k-2},I_{k-2})$... Continuing this 
procedure inductively we finally get the desired elements $h_i$ 
and intervals $I_j$, thus concluding the proof of the lemma.

\vspace{0.3cm}

From now on we assume that $\Gamma$ is not solvable with degree of 
solvability less than or equal to $k$. Note that in this case the 
elements $h_i$ constructed in the preceding lemma do not necessarily 
satisfy a property of ``periodicity'' so strong as (\ref{call}) 
or (\ref{compli}). Nevertheless, as the next lemma shows, 
they must be ``non lacunary'' in a very precise sense.

\vspace{0.1cm}

\begin{lem} {\em For every $i \in \{2,\ldots,k\}$ the set of integers}
$$\mathcal{N}_i = \{ m \geq 0: \quad \!\!\! h_i 
\mbox{ has no fixed point inside } h_1^m (I_{i-1}) \}$$
{\em is syndetic, i.e. it cannot have arbitrarily large gaps.}
\end{lem}

\noindent{\bf Proof.} Assuming that $\mathcal{N}_i$ has arbitrarily large 
gaps, we will prove that $h_1$ and $h \!=\! h_i$ generate a free 
semi-group on two generators, thus giving a contradiction. For 
this let us consider two words \esp \esp 
$W_1 = h_1^{n} h^{m_r} h_1^{n_r} \cdots h^{m_1} h_1^{n_1}$ \esp \esp 
and \esp \esp $W_2 = h^{q} h_1^{p_s} h^{q_s} \cdots h_1^{p_1} h^{q_1}$ 
in positive powers of $h_1$ and $h$, \esp 
\esp where $m_j,n_j,p_j,q_j$ are positive integers, $n \geq 0$, and 
$q \geq 0$ (with $n > 0$ when $r=0$, and $q > 0$ when $s=0$). We have 
to prove that these words represent different elements of $\Gamma$.

First of all, using a similar argument to that of the beginning of the 
proof of Lemma \ref{max}, one easily checks that if the numbers 
\esp $N_1 = n_1 + \ldots + n_r + n$ \esp and $N_2 = p_1 + \ldots + p_s$ 
are different, then $W_1$ and $W_2$ are distinct elements. Assume in what 
follows that $N_1$ and $N_2$ coincide with some $N \in \mathbb{N}$. We will 
finish the proof by checking that $h_1^{-N} W_1$ and $h_1^{-N} W_2$ are 
different. For this note that
\begin{eqnarray*}
h_1^{-N} W_1 
&=& h_1^{-N} h_1^{n} h_i^{m_r} h_1^{n_r} \cdots h_i^{m_1} h_1^{n_1} \\
&=& h_1^{-N} h_1^{n} h_i^{m_r} h_1^{n_r} \cdots h_i^{m_2} h_1^{n_1+n_2} 
\big( h_1^{-n_1} h_i^{m_1} h_1^{n_1} \big) \\
&=& h_1^{-N} h_1^{n} h_i^{m_r} h_1^{n_r} \cdots h_i^{m_3} h_1^{n_1+n_2+n_3} 
\big( h_1^{-(n_1+n_2)} h_i^{m_2} h_1^{n_1+n_2} \big) \esp 
\big( h_1^{-n_1} h_i^{m_1} h_1^{n_1} \big) \\
&\vdots& \\
&=& \big( h_1^{-(N-n)} h_i^{m_r} h_1^{N-n} \big) \cdots  
\big( h_1^{-(n_1+n_2)} h_i^{m_2} h_1^{n_1+n_2} \big) \esp 
\big( h_1^{-n_1} h_i^{m_1} h_1^{n_1} \big),
\end{eqnarray*} 
and 
\begin{eqnarray*}
h_1^{-N} W_2 
&=& h_1^{-N} h_i^{q} h_1^{p_s} h^{q_s} \cdots h_1^{p_2}h_i^{q_2}h_1^{p_1}h_i^{q_1} \\
&=& h_1^{-N} h_i^{q} h_1^{p_s} h_i^{q_s} \cdots 
h_i^{q_3} h_1^{p_1+p_2} \big( h_1^{-p_1} h_i^{q_2} h_1^{p_1} \big) h_i^{q_1} \\
&\vdots& \\
&=& \big( h_1^{-N} h_i^q h_1^{N} \big) \esp 
\big( h_1^{-(N - p_s)} h_i^{q_s} h_1^{N - p_s} \big) \cdots 
\big( h_1^{-p_1} h_i^{q_2} h_1^{p_1} \big) h_i^{q_1}.
\end{eqnarray*}
Recall that, by the definition of the star operation, for every 
$n \geq 0$ the interval $h_1^n(I_{i-1})$ is fixed by $h_i$. Now take 
$m \in \mathcal{N}_i$ such that the element next to it in $\mathcal{N}_i$ 
is bigger than $m + N$. For this choice of $m$, all the elements 
$$\big( h_1^{-(N-n)} h_i^{m_r} h_1^{N-n} \big), \esp 
\big( h_1^{-(n_1+n_2)} h_i^{m_2} h_1^{n_1+n_2} \big), \esp 
\ldots, \esp \big( h_1^{-n_1} h_i^{m_1} h_1^{n_1} \big)$$ 
and
$$\big( h_1^{-N} h_i^q h_1^{N} \big), \esp  
\big( h_1^{-(N - p_s)} h_i^{q_s} h_1^{N - p_s} \big), \esp 
\ldots, \esp 
\big( h_1^{-p_1} h_i^{q_2} h_1^{p_1} \big), \esp h_i^{q_1}$$
fix the interval $h_1^m(I_{i-1})$, but only the last one (namely 
$h_i^{q_1}$) has no fixed point inside it. Therefore, if we let $\nu$ 
be any Radon measure on this interval which is invariant by the group 
generated by (the restrictions) of all these maps, then one has 
$$\tau_{\nu}(h_1^{-N} W_1) = 0 \quad \mbox{ and } \quad 
\tau_{\nu}(h_1^{-N} W_2) = \tau_{\nu}(h_i^{q_1}) \neq 0.$$
This shows that \esp $h_1^{-N} W_1 \neq h_1^{-N}W_2$ 
\esp and finishes the proof of the lemma.

%%%%%%%%%%%%%%%%%%%%%%%%%%%%%%%%%%%%%%%%%%%%%%%%%%%%%%%%%%%%%%%%%%%%%%%%%%%%%%%%%
%%%%%%%%%%%%%%%%%%%%%%%%%%%%%%%%%%%%%%%%%%%%%%%%%%%%%%%%%%%%%%%%%%%%%%%%%%%%%%%%%

\subsubsection{End of the proof}
\label{fin}

\hspace{0.35cm} We are now 
ready to finish the proof of Theorem B. For this first note 
that the preceding lemma shows the existence of a (large) positive integer $N$ 
satisfying the following property: for every $i \in \{2,\ldots,k\}$ and every 
integer $m \geq 0$, at least one of the maps  
$$h_i, h_1^{-1} h_i h_1^{1}, \ldots, h_1^{-N} h_i h_1^{N}$$ 
has no fixed point on the interval $h_1^m (I_{i-1})$.

Now proceed as in the proof of Proposition \ref{cagon}. Let $M$ be 
a common $\alpha$-H\"older constant for the function $\log(h_1')$ 
and for all the maps of the form \esp $ h_1^{-j} h_i h_1^j$, \esp where 
\esp $j \in \{0,1,\ldots,N\}$ \esp and \esp $i \in \{2,\ldots,k-2\}$. 
\esp Let $\bar{M}$ now be defined by
$$\log(\bar{M}) = 
\big( 1 + \alpha [1 + (1+\alpha) + (1+\alpha)^2 + \ldots + (1+\alpha)^{k-4}] \big) 
\esp \log(M).$$ 
For simplicity let us denote $]a_i,b_i[ = I_{i}$ for every 
$i \in \{1,\ldots,k-1\}$, and let us fix a constant $\lambda$ such that 
$$1 < \lambda < 1 + \frac{b_{k-1} - a_{k-1}}{e^{\bar{M}} \esp (a_{k-1} - a_{k-2})}.$$
The derivative of $h_{k-1}$ must be equal to $1$ at the origin, which is the 
accumulation point of the intervals $h_1^m([a_{k-2},b_{k-2}])$ for $m \geq 0$. 
Therefore, we can fix a (large) positive integer $n$ such that 
$(h_1^{-j} h_{k-1} h_1^j)' (x) \leq \lambda$ and 
$(h_1^{-j} h_{k-1}^{-1} h_1^j)'(x) \leq \lambda$ for all 
$x \in h_1^n ([a_{k-2},b_{k-2}])$ and all $j \in \{ 0 ,1 ,\ldots, N \}$.
Since $]a_{k-1},b_{k-1}[$ belongs to $\mathcal{I}(h_{k})$ 
and $\Gamma$ has no crossed elements, one of the maps 
$$h_{k-1},h_1^{-1}h_{k-1}h_1,\ldots,h_1^{-N}h_{k-1}h_1^N$$
must send the interval $h_1^n(]a_{k-1},b_{k-1}[)$ 
into a disjoint interval (still contained in 
$h_1^n(]a_{k-2},b_{k-2}[)$). Assume for instance that this element 
is $h_{k-1}$ (all the other cases are 
treated similarly), 
and that $h_{k-1}(h_1^n(]a_{k-1},b_{k-1}[))$ is to 
the right of $h_1^n(]a_{k-1},b_{k-1}[)$, {\em i.e.} 
that is $h_{k-1}(h_1^n(a_{k-1})) \geq h_1^n(b_{k-1})$ 
(if not then replace $h_{k-1}$ by its inverse). 
There must exist points 
$u \!\in\! [a_{k-1},b_{k-1}] \subset [a_{k-2},b_{k-2}]$ 
and $v \!\in\! [a_{k-2},a_{k-1}] \subset [a_{k-2},b_{k-2}]$ such that
\begin{equation}
\frac{h_{k-1}(h_1^{n}(a_{k-1})) - h_1^{n}(a_{k-1})}{h_1^{n}(a_{k-1}) -
h_1^{n}(a_{k-2})} 
\geq \frac{h_1^{n}(b_{k-1}) - h_1^{n}(a_{k-1})}{h_1^{n}(a_{k-1}) - h_1^{n}(a_{k-2})} = 
\frac{(h_1^{n})'(u)}{(h_1^{n})'(v)} \cdot \frac{b_{k-1}-a_{k-1}}{a_{k-1} - a_{k-2}}.
\label{mani2}
\end{equation}
Note that 
\begin{eqnarray*}
\left| \log \Big( \frac{(h_1^{n})'(u)}{(h_1^{n})'(v)} \Big) \right| 
&=& \left| \log \Big( \prod_{m=0}^{n - 1} \frac{h_1'(h_1^m(u))}{h_1'(h_1^m(v))}
\Big) \right| \\
&\leq& \sum_{m=0}^{n - 1} \big| \log(h_1'(h_1^m(u))) - \log(h_1'(h_1^m(v))) \big| \\
&\leq& M \esp \sum_{m=0}^{n - 1} \big| h_1^m(u) - h_1^m(v) \big|^{\alpha} \\
&\leq& M \esp \sum_{m=0}^{n - 1} \big| h_1^m(b_{k-2}) - h_1^m (a_{k-2})
\big|^{\alpha}. 
\end{eqnarray*}
Now for each $m \!\in\! \{0,1,\ldots,n-1\}$ there exists $j \!\in\! \{0,1,\ldots,N\}$ 
(depending on $m$) such that the intervals $h_1^m \big( ]a_{k-2},b_{k-2}[ \big)$ 
and $h_1^{-j} h_{k-2} h_1^j \big( h_1^m(]a_{k-2},b_{k-2}[) \big)$ are disjoint. 
Lemma \ref{cagoncito} applied to \esp $h_1^{-j} h_{k-2} h_1^j$ \esp gives
$$|h_1^m(b_{k-2}) - h_1^m(a_{k-2})| 
\leq M \esp |h_1^m(b_{k-3}) - h_1^m(a_{k-3})|^{1+\alpha},$$
and so
$$\left| \log \Big( \frac{(h_1^{n})'(u)}{(h_1^{n})'(v)} \Big) \right| \leq 
M \esp \sum_{m=0}^{n-1} M^{\alpha} \esp 
|h_1^m(b_{k-3}) - h_1^m(a_{k-3})|^{\alpha(1+\alpha)}.$$ 
Repeating this argument several times we deduce
\begin{eqnarray*}
\left| \log \Big( \frac{(h_1^{n})'(u)}{(h_1^{n})'(v)} \Big) \right| 
&\leq& M^{1+\alpha} \sum_{m=0}^{n-1} M^{\alpha(1+\alpha)} 
\big| h_1^m(b_{k-4}) - h_1^m(a_{k-4})\big|^{\alpha(1+\alpha)^{2}}\\
&\vdots& \\
&\leq& M^{1+\alpha+\alpha[(1+\alpha)+\ldots+(1+\alpha)^{k-4}]} \esp 
\sum_{m=0}^{n-1} \big| h_1^m(b_1) - h_1^m(a_1)\big|^{\alpha(1+\alpha)^{k-3}},
\end{eqnarray*}
that is
\begin{equation}
\left| \log \Big( \frac{(h_1^{n})'(u)}{(h_1^{n})'(v)} \Big) \right| \leq 
\bar{M} \esp 
\sum_{m=0}^{n-1} \big| h_1^m(b_1)-h_1^m(a_1) \big|^{\alpha(1+\alpha)^{k-3}}.
\label{chita2}
\end{equation}
By the definition of $k$ (namely $k = 1+k(\alpha)$), 
one has $\alpha(1+\alpha)^{k-3} \geq 1$. Moreover, the 
intervals $h_1^m (]a_{1},b_{1}[)$ must be pairwise disjoint. 
Hence, the right hand side expression of (\ref{chita2}) is bounded by 
$\bar{M}$, which implies that 
$$\frac{(h_1^{n})'(u)}{(h_1^{n})'(v)} \geq \frac{1}{e^{\bar{M}}}.$$
Introducing this inequality into (\ref{mani2}) we get 
$$\frac{h_k(h_1^{n}(a_{k-1})) - h_1^{n}(a_{k-1})}{h_1^{n}(a_{k-1}) - h_1^{n}(a_{k-2})} 
\geq \frac{b_{k-1} - a_{k-1}}{e^{\bar{M}} \esp (a_{k-1} - a_{k-2})},$$
and summing 1 to both members this gives
$$\frac{h_{k-1}(h_1^{n}(a_{k-1}))-h_1^{n}(a_{k-2})}{h_1^{n}(a_{k-1})-h_1^{n}(a_{k-2})} 
\geq 1 + \frac{b_{k-1} - a_{k-1}}{e^{\bar{M}} \esp (a_{k-1} - a_{k-2})} > \lambda.$$
However, since $h_{k-1}$ fixes the point $h_1^n(a_{k-2})$, the left hand 
side member of this inequality is equal to $h_{k-1}'(x)$ for some point 
$x \in [h_1^n(a_{k-2}),h_1^n(a_{k-1})]$, and so it is less than or equal 
to $\lambda$ by our choice of $n$. This contradiction 
finishes the proof of Theorem B.

\vspace{0.1cm}

\begin{rem} A careful reading of the arguments given along the second part 
of this article shows that the differentiability of the maps involved 
is needed only at one of the end points of $[0,1]$. More precisely, 
Theorem B still holds (with the very same proof) for finitely 
generated subgroups of $\mathrm{Diff}_+^{1+\alpha}([0,1[)$ or 
$\mathrm{Diff}_+^{1+\alpha}(]0,1])$ without free semi-groups 
on two generators. This also applies to Remark \ref{clog}. 
However, we ignore if the theorem is still true for 
groups of germs of $C^{1+\alpha}$ diffeomorphisms; 
this seems to be an interesting problem.
\label{un-solo-lado}
\end{rem}

%%%%%%%%%%%%%%%%%%%%%%%%%%%%%%%%%%%%%%%%%%%%%%%%%%%%%%%%%%%%%%%%%%%%%%%%%%%%%%%%%%%%%%%%%%%%%%%
%%%%%%%%%%%%%%%%%%%%%%%%%%%%%%%%%%%%%%%%%%%%%%%%%%%%%%%%%%%%%%%%%%%%%%%%%%%%%%%%%%%%%%%%%%%%%%%

\subsection{The cases of the circle and the real line}
\label{ext}

\hspace{0.35cm} The aim of this final section is to prove two claims made 
in the Introduction of this work, namely that finitely generated 
subgroups of $\mathrm{Diff}_+^{1+\alpha}(\mathbb{R})$ or 
$\mathrm{Diff}^{1+\alpha}_+(\mathrm{S}^1)$ with 
sub-exponential growth are also almost nilpotent. 
Again, we will prove this for subgroups without free semi-groups 
on two generators by showing that they are solvable with degree of 
solvability at most $2 + k(\alpha)$. (Remark that this issue will 
still be true for non finitely generated groups without free 
semi-groups on two generators.) The nilpotence will then be 
a direct consequence of Rosenblatt's theorem \cite{Ros}.

Let us first consider the (simpler) case of a finitely generated subgroup 
$\Gamma$ of $\mathrm{Diff}_+^{1+\alpha}(\mathbb{R})$ without free semi-groups 
on two generators. By \S \ref{bp}, the action of the first derived group 
$\Gamma_1$ has global fixed points. Looking at the action of 
$\Gamma_1$ on the closure of each connected component of the 
complement of the set of its global fixed points, and using (the 
arguments of the proof of) Theorem B (as well as Remark 
\ref{un-solo-lado}), we obtain that $\Gamma_1$ is solvable 
with degree of solvability at most $1 \! + \! k(\alpha)$. 
Then one deduces that $\Gamma$ is solvable itself with degree 
of solvability smaller than or equal to $2 \! + \! k(\alpha)$.

Now let $\Gamma$ be a finitely generated subgroup of 
$\mathrm{Diff}_+^{1+\alpha}(\mathrm{S}^1)$ without 
free semi-groups on two generators. Note that, 
in contrast to the case of sub-exponential growth 
subgroups, $\Gamma$ is not {\em a priori} amenable.\footnote{There 
exist non amenable groups without free semi-groups on two generators 
(see for example \cite{olsha}), 
but it seems to be unknown if such a group can act faithfully on 
the circle.} This is the reason why the following Claim is not 
completely trivial.\footnote{According to a beautiful result 
by Margulis \cite{mar}, the Claim still holds when $\Gamma$ 
has no free subgroup on two generators.}

\vspace{0.32cm}

\noindent{\bf{\underbar{Claim:}}} The group $\Gamma$ 
preserves a probability measure on the circle.

\vspace{0.2cm}

\noindent{\bf Proof.} Let $\tilde{\Gamma}$ be the covering of $\Gamma$ 
acting on the real line. This group $\tilde{\Gamma}$ is still finitely 
generated. Moreover, it cannot contain crossed elements: if it contains 
two such elements then they project on $\Gamma$ into two elements 
for which one can apply the argument of Proposition \ref{semilibre} 
in order to show that $\Gamma$ contains free semi-groups 
on two generators, thus giving a contradiction. 

By Proposition \ref{radon}, the group $\tilde{\Gamma}$ preserves a 
Radon measure on the real line. This measure is invariant by the integer 
translations, and so it projects into a finite measure on the circle 
which is invariant by $\Gamma$. Hence, up to normalization, we 
have obtained the desired $\Gamma$-invariant probability 
measure on $\mathrm{S}^1$.

\vspace{0.28cm}

For groups of circle homeomorphisms preserving a probability measure 
it is easy to see that the {\em rotation number} function is a 
homomorphism into $\mathbb{T}^1$. (See for instance \cite{ghys}.) 
Using this fact one easily deduces that in our situation the first 
derived group $\Gamma_1$ has global fixed points. Again, looking at 
the action of $\Gamma_1$ on the closure of each connected component 
of the complement of the set of its global fixed points, and using (the 
arguments of the proof of) Theorem B, we obtain that $\Gamma_1$ is 
solvable with degree of solvability at most $1 \! + \! k(\alpha)$. 
From this one concludes that $\Gamma$ is solvable itself with 
degree of solvability less than or equal to $2 \! + \! k(\alpha)$, 
thus finishing the proof.

%%%%%%%%%%%%%%%%%%%%%%%%%%%%%%%%%%%%%%%%%%%%%%%%%%%%%%%%%%%%%%%%%%%%%%%%%%%%%%%%%%%%%%%%%%%%%%%
%%%%%%%%%%%%%%%%%%%%%%%%%%%%%%%%%%%%%%%%%%%%%%%%%%%%%%%%%%%%%%%%%%%%%%%%%%%%%%%%%%%%%%%%%%%%%%%

\begin{small}

\vspace{0.1cm}

\noindent Andr\'es Navas\\

\noindent Univ. de Santiago de Chile, Alameda 3363, Santiago, Chile\\

\noindent Univ. de Chile, Las Palmeras 3425, \~Nu\~noa, Santiago, Chile\\

\end{small}

\end{document}